\documentclass{article}
\usepackage{amsmath}
\usepackage{amssymb}
\usepackage{amsthm}
\usepackage{amsfonts}
\usepackage{enumitem}
\usepackage[dvipsnames]{xcolor}
\usepackage{hyperref}
\usepackage[capitalise]{cleveref}
\usepackage{graphicx} 
\usepackage{tikz-cd}
\usepackage{quiver}
\usepackage{tikzit}
\usepackage{adjustbox}

\DeclareFontFamily{U}{mathx}{}
\DeclareFontShape{U}{mathx}{m}{n}{<-> mathx10}{}
\DeclareSymbolFont{mathx}{U}{mathx}{m}{n}
\DeclareMathAccent{\widehat}{0}{mathx}{"70}
\DeclareMathAccent{\widecheck}{0}{mathx}{"71}
\DeclareMathAccent{\widebar}{0}{mathx}{"73} 

\newcommand{\R}{\mathbb{R}}
\newcommand{\g}{\mathfrak{g}}
\newcommand{\ad}{\mathrm{ad}}
\newcommand{\Ad}{\mathrm{Ad}}
\newcommand{\Sym}{\mathrm{Sym}}
\DeclareMathOperator{\coker}{coker}
\DeclareMathOperator{\rk}{rk}
\newcommand{\VF}{\mathfrak{X}}
\DeclareMathOperator{\hofib}{hofib}
\DeclareMathOperator{\clim}{colim}
\DeclareMathOperator{\colim}{hocolim}
\DeclareMathOperator{\holim}{holim}
\DeclareMathOperator{\Der}{Der}
\newcommand{\cdga}{\mathsf{cdga}}
\newcommand{\Mod}{\mathsf{\text{-}Mod}}

\newcommand{\Tang}{\mathbb{T}}
\newcommand{\Cotang}{\mathbb{L}}

\newcommand{\Spec}{\mathrm{Spec}}
\newcommand{\Map}{\mathrm{Map}}

\newcommand{\DecompCircle}{\begin{tikzpicture}[none/.style={}, dot/.style={fill=black, draw=black, shape=circle, scale=0.5}, white-dot/.style={fill=white, draw=black, shape=circle, scale=0.5}, to/.style={->}, scale=0.3, baseline=.1ex]
	\begin{pgfonlayer}{nodelayer}
		\node [style=dot] (0) at (0, 1) {};
		\node [style=white-dot] (1) at (0, -1) {};
		\node [style=dot] (2) at (5, 3) {};
		\node [style=white-dot] (3) at (5, 1) {};
		\node [style=dot] (4) at (4, -1) {};
		\node [style=white-dot] (5) at (4, -3) {};
		\node [style=none] (6) at (1, -0.5) {};
		\node [style=none] (7) at (3.25, -1.5) {};
		\node [style=none] (8) at (1, 0.5) {};
		\node [style=none] (9) at (3.25, 1.5) {};
	\end{pgfonlayer}
	\begin{pgfonlayer}{edgelayer}
		\draw [bend left=90, looseness=1.75] (4) to (5);
		\draw [bend right=270, looseness=1.75] (3) to (2);
		\draw [style=to] (8.center) to (9.center);
		\draw [style=to] (6.center) to (7.center);
	\end{pgfonlayer}
\end{tikzpicture}}

\newcommand{\BettiCircle}{\begin{tikzpicture}[none/.style={}, dot/.style={fill=black, draw=black, shape=circle, scale=0.5}, white-dot/.style={fill=white, draw=black, shape=circle, scale=0.5}, to/.style={->}, scale=0.3, baseline=-0.5ex]
	\begin{pgfonlayer}{nodelayer}
		\node [style=dot] (0) at (0, 1) {};
		\node [style=dot] (1) at (0, -1) {};
		\node [style=dot] (4) at (4.5, -2) {};
		\node [style=none] (6) at (1, -0.5) {};
		\node [style=none] (7) at (3.25, -1.5) {};
		\node [style=none] (8) at (1, 0.5) {};
		\node [style=none] (9) at (3.25, 1.5) {};
		\node [style=dot] (10) at (4.5, 2) {};
        \node [style=none] (11) at (5, 0) {};
	\end{pgfonlayer}
	\begin{pgfonlayer}{edgelayer}
		\draw [style=to] (8.center) to (9.center);
		\draw [style=to] (6.center) to (7.center);
	\end{pgfonlayer}
\end{tikzpicture}
}

\newcommand{\BettiCircleb}{\begin{tikzpicture}[none/.style={}, dot/.style={fill=black, draw=black, shape=circle, scale=0.5}, white-dot/.style={fill=white, draw=black, shape=circle, scale=0.5}, to/.style={->}, scale=0.3, baseline=-0.5ex]
	\begin{pgfonlayer}{nodelayer}
		\node (0) at (0, 1) {*};
		\node (1) at (0, -1) {*};
		\node (4) at (4.5, -2) {*};
		\node [style=none] (6) at (1, -0.5) {};
		\node [style=none] (7) at (3.25, -1.5) {};
		\node [style=none] (8) at (1, 0.5) {};
		\node [style=none] (9) at (3.25, 1.5) {};
		\node (10) at (4.5, 2) {*};
        \node [style=none] (11) at (5, 0) {};
	\end{pgfonlayer}
	\begin{pgfonlayer}{edgelayer}
		\draw [style=to] (8.center) to (9.center);
		\draw [style=to] (6.center) to (7.center);
	\end{pgfonlayer}
\end{tikzpicture}
}

\title{Shifted Symplectic Geometry by Examples}
\author{Damien Calaque and Stefano Ronchi}
\date{}

\newtheorem{Theorem}{Theorem}[section]
\newtheorem{Idea}[Theorem]{Idea}

\theoremstyle{definition}
\newtheorem{Remark}[Theorem]{Remark}
\newtheorem{Definition}[Theorem]{Definition}
\newtheorem{Lemma}[Theorem]{Lemma}
\newtheorem{Example}[Theorem]{Example}
\newtheorem{Exercise}[Theorem]{Exercise}
\newtheorem{Notation}[Theorem]{Notation}

\begin{document}

\maketitle

\setcounter{tocdepth}{2}
\tableofcontents


\section{Introduction}

These notes are intended to be an introduction to shifted symplectic geometry, targeted to Poisson geometers with a serious background in homological algebra. They are extracted from a mini-course given by the first author at the Poisson 2024 summer school that took place at the Accademia Pontaniana in Napoli. 

\medskip

It is worth noticing that shifted symplectic geometry (with a non-negative shift) was known to Poisson geometers before it even existed. For instance, Xu's quasi-symplectic groupoids \cite{Xu} are exactly $1$-shifted symplectic underived smooth $1$-stacks in the differentiable setting, as we will see in these notes. Similarly, twisted Courant algebroids and twisted Dirac structures can be understood as $2$-shifted symplectic formal stacks and lagrangian morphisms into those (see \cite{Pym-Safronov} and references therein). 

Instances of negatively shifted symplectic structures have also been explored by that part of the mathematical physics community that is very familiar to (and with) Poisson geometers. For example, the ``antibracket'' of the BV-BRST formalism is a by-product of the $(-1)$-shifted symplectic structure on derived critical loci. The ``AKSZ'' transgression procedure (named after \cite{AKSZ}) led to the striking result (of \cite{PTVV}, where general shifted symplectic structures were first introduced) about the existence of shifted symplectic structures on mapping stacks. The AKSZ construction being relevant to both classical and quantum field theory, it is not surprising that shifted symplectic structures also are related to classical field theories (see e.g.~\cite{Cal,CalSurvey,CHS}) and quantization (see \cite{PTVV,ToenSurvey,CPTVV}). 

Still in connection with mathematical physics, algebraic geometers working in enumerative geometry were also very close to discovering (negatively) shifted symplectic structures. For instance, the symmetric obstruction theories of Behrend--Fantechi \cite{BF2} are the shadow of $(-1)$-shifted symplectic structures, and their work on lagrangian intersections \cite{BF1} witnesses another striking result (of \cite{PTVV} again) saying that derived lagrangian intersections are shifted symplectic. 
We refer to \cite{Joyce} (especially Figure 1.1) for more details on the role of $(-1)$-shifted symplectic structures in Donaldson--Thomas theory. 

\medskip

In order to have a global, intrinsic, and model-independent approach to shifted symplectic structures, it seems unavoidable to use the abstract language of derived geometry (see e.g.~\cite{HAG-II}). This approach is taken in the aforementioned seminal paper \cite{PTVV} by Pantev--Toën--Vaquié--Vezzosi, which introduced shifted symplectic structures on derived stacks, and it allows them to prove very general structural results. It has several crucial advantages: 
\begin{itemize}
\item It makes it easier to prove general results, and even to \textit{define} things in full generality. 
\item It is intrinsic: geometric objects are often defined in a very natural way, and more traditional constructions are obtained by \textit{computation}. It then becomes tautological that different constructions/computations lead to the same object. 
\item It works under very mild assumptions, and thus it allows dealing with rather singular spaces. 
\end{itemize}
But everything has a cost: 
\begin{itemize}
    \item The foundations of derived geometry use a lot of abstract homotopy theory (including $\infty$-categories in more modern expositions). This means that there is a lot of material that one has to learn to understand even the most basic definitions. 
    \item Going back to examples and explicit computations within a specific model can be difficult\footnote{Note that the difficulties with explicit computations are often not
a purely derived problem, but appear equally well at the classical level: objects with a clean moduli space description often have a local geometry that is poorly understood.}. 
\end{itemize}
As a result, the literature on shifted symplectic geometry is often perceived as too abstract and not so much connected to the interests of many geometers. 

Despite this, there has recently been some really great work in Poisson geometry that makes use of the language and intuition of shifted symplectic geometry (see e.g.~\cite{Cueca-Zhu,Crooks-Mayrand} and references therein), while still technically relying on ``good old'' models. It seems that the gap between communities remains. This set of lecture notes does not at all pretend to fill this gap, but tries to provide a bottom-up introduction to (the abstract approach to) shifted symplectic geometry. 

We particularly provide computations that are not made explicit in most references (as the intrinsically homotopical formalism is somehow used to deal with these computations ``by itself, in the background''). Along the way, we hope to convince the reader that (a) derived stacks are unavoidable at some point, (b) abstract definitions can be made very concrete in examples, (c) there is still a lot of interesting comparison work that remains to be done even in well-known examples\footnote{On the quantum side of the story, a similar type of comparison work has been carried out with great success by Ben-Zvi--Borchier--Jordan \cite{BZBJ}, who compared previously known explicit quantizations of character varieties to the quantization obtained using factorization homology (which can be seen as a non-commutative analog of the mapping stack construction). It is not a coincidence that shifted symplectic geometry and factorization homology have emerged around the same years, and that abstract homotopy theory serves as foundations for both. }.

\medskip

We would finally like to advertise another set of lecture notes (both recent and inspiring) by Cueca--Maglio--Valencia \cite{CMV}, that will provide a perfect complement to the present ones. 

\subsection*{Organization of the paper}

\cref{section2} provides an introduction to shifted symplectic structures and lagrangian structures in the linear setting\footnote{Considerably expending the discussion from the survey \cite{CalSurvey}.}. We apply the definition to fairly simple examples of cochain complexes, such as $1$-term and $2$-term complexes, and try to make all homotopies as explicit as we can. We also introduce some specific features and examples that will appear more systematically in the geometric setting. 

\cref{section3} introduces shifted symplectic structures on derived affine schemes and lagrangian structures on morphisms thereof. Before doing so, we explain the main features of the homotopy theory of commutative differential graded algebras that is required for the reader to understand what is going on: (co)tangent complexes, derived tensor products (which geometrically correspond to derived intersections), etc...

In principle, \cref{section4} deals with shifted symplectic structures on \textit{derived higher} stacks, but a large part of the discussion (and examples) is about underived $1$-stacks. We provide an exhaustive description of shifted symplectic structures on such stacks that are presented as quotient of Lie groupoids, and relate these to notions that are familiar to Poisson geometers. We then proceed similarly for lagrangian morphisms, and explain how they are incarnations of various notions of moment maps. We finally interpret the reduction procedure for these moment maps in terms of lagrangian intersections, and see our first examples of derived stacks appearing. 

Finally, we present in \cref{section5} the so-called AKSZ construction for shifted symplectic structures (also called PTVV after \cite{PTVV}). It allows to construct shifted symplectic structures on (derived) mapping stacks. This section is more sketchy and provides less detailed computations (in order to keep these notes both readable and within a reasonable length). It should be understood as an invitation to more advanced topics in shifted symplectic geometry, and to get one's hands on comparison questions between the new abstract constructions and the more concrete ones. We use character varieties/stacks as a leading example. 

\subsection*{Conventions}
All along these notes, $k$ is a field of characteristic zero. We use $\cong$ for isomorphisms, and $\simeq$ for quasi-isomorphisms (and more generally for all kinds of equivalences that are weaker than genuine isomorphisms: Morita equivalences, homotopy equivalences, equivalences in an $\infty$-category, etc...). 

\subsection*{Acknowledgments}
DC warmly thanks the organizers of the Poisson 2024 summer school for giving him the opportunity to deliver this series of lectures, and more generally for a fantastic and very successful week. Both authors thank the participants for their enthusiasm. 
DC also thanks Tristan Bozec for his careful reading and for his comments on a preliminary version of these notes. Finally, both authors thank the anonymous referee for their excellent suggestions. 


\section{Shifted symplectic linear algebra}\label{section2}

\subsection{Main ideas}

Recall that a \emph{symplectic structure} on a $k$-vector space $V$ is a linear map $\omega:\wedge^2V\to k$ that is \emph{non-degenerate}, meaning that  
\begin{eqnarray*}
\omega^\flat:V  & \longrightarrow   & V^* \\
v               & \longmapsto       & \omega(v\wedge-)
\end{eqnarray*}
is an isomorphism. 

\begin{Idea}\label{first idea}
Replace $k$-vector spaces with cochain complexes (of $k$-vector spaces) and isomorphisms with quasi-isomorphisms. 
\end{Idea}
Let $(V,\omega)$ a symplectic $k$-vector space and $L\subset V$ a vector subspace. The subspace $L$ is said \emph{lagrangian} if $\omega|_{\wedge^2L}=0$ (in which case one says that $L$ is \emph{isotropic}) and $L$ is maximal for this property. The maximality property can be equivalently rephrased in the following (equivalent) ways: 
\begin{itemize}
\item[(a)] $\dim_k(L)=\frac12\dim_k(V)$; 
\item[(b)] The inclusion $L\subset L^o:=\{v\in V\,|\,\omega^\flat(v)_{|L}=0\}$ is an equality. 
\item[(c)] The null-sequence\footnote{A ``null-sequence'' is just another name for a complex: composing two maps in the sequence gives zero. } $0\to L\to V\simeq V^*\to L^*\to 0$ is exact. 
\end{itemize}

How is \cref{first idea} incarnated in this situation?
\begin{Idea}\label{second idea}
We will replace $\omega_{|\wedge^2L}=0$ by the condition that $\omega_{|\wedge^2L}:\wedge^2L\to k$ is \emph{homotopic} to zero. 
The data of the homotopy $h$ will be part of the structure.  
This is equivalent to requiring that the composition 
\begin{equation}\label{null-sequence}
(\omega_{|\wedge^2L})^\flat:L\longrightarrow V\overset{\omega^\flat}{\longrightarrow} V^*\longrightarrow L^*
\end{equation}
is homotopic to zero (we will say \emph{null-homotopic}), via the homotopy $h^\flat$.  
\end{Idea}
\begin{Idea}\label{third idea}
Recall from \cref{first idea} that we assume $\omega^\flat$ is a quasi-isomorphism. 
Inspired by its equivalent formulation (c), the maximality condition shall be replaced with the condition that the null-homotopic sequence \eqref{null-sequence} is an exact triangle. In particular, it induces a long exact sequence
\begin{equation}\label{long sequence}
\cdots\to H^\bullet(L)\to H^\bullet(V)\simeq H^\bullet(V^*)\to H^\bullet(L^*)\to H^{\bullet+1}(L)\to\cdots
\end{equation}
in cohomology. 
\end{Idea}

\begin{Remark}\label{remark-context}
All the above still makes sense in an arbitrary symmetric monoidal abelian $k$-linear category, such as that of $A$-modules ($A$ being a commutative $k$-algebra), that of $G$-representations ($G$ being a group), ... 

Note that, even though the symmetric monoidal $\mathbb{R}$-linear category of vector bundles over a smooth manifold $M$ is not abelian, the category of complexes thereof is good enough for our purposes. 

The most suitable general context for what we do here is the one of symmetric monoidal stable $k$-linear $\infty$-categories; but we will avoid the language of $\infty$-categories as much as we can. 
\end{Remark}

\subsection{Linear shifted symplectic structures}\label{section:linear}

Let $V$ be a complex of $k$-vector spaces (or $A$-modules, $G$-representations, etc... see \cref{remark-context}). 

\begin{Definition}
    An \emph{$n$-shifted symplectic structure} on $V$ is a cochain map $\omega: \wedge^2 V \to k[n]$ such that
    \begin{eqnarray*}
    \omega^\flat:V  & \longrightarrow   & V^*[n] \\
    v               & \longmapsto       & \omega(v\wedge-)
    \end{eqnarray*}
    is a quasi-isomorphism. The condition that $\omega^\flat$ is a quasi-isomorphism is called \emph{non-degeneracy condition}. 
\end{Definition}

Recall that $V[n]$ is the cochain complex $V$ shifted by $n$, that is $V[n]^k = V^{k+n}$ (and its differential is $(-1)^n$ times the differential of $V$). If $V$ is concentrated in degree $0$, then $V[n]$ is concentrated in degree $-n$. For complexes, the second exterior power is defined as $\wedge^2V:=\mathrm{Sym}^2(V[-1])[2]$, where $\mathrm{Sym}^2(W)$ is the quotient of $W^{\otimes 2}$ by the relation $a\otimes b=(-1)^{|a||b|}b\otimes a$. 

\begin{Remark}
    In general $k$ is replaced by the monoidal unit of the category where the cochain complexes are based. This is for example the algebra $A$ for $A$-modules, the trivial character for representations of a group $G$, \dots
\end{Remark}

\begin{Remark}
    In general, to be in line with \cref{first idea}, we should be considering forms $\wedge^2 \widetilde{V} \to k[n]$, for any complex $\widetilde{V}$ quasi-isomorphic to $V$. 
    But if $V$ is nice enough (e.g. made of projectives) this is equivalent. 
\end{Remark}

\begin{Example}\label{ex:shift-sympl-struct-deRham-cplx}
    Let $X$ be an $n$-dimensional closed oriented manifold. Consider the cochain complex of differential forms with the de Rham differential shifted by 1: $V=(\Omega^\bullet(X), d_{dR})[1]$. 
    This is $(2-n)$-shifted symplectic with respect to the form 
    \begin{equation*}
        \omega (\alpha \wedge \beta) := \int_X \alpha \wedge \beta,
    \end{equation*}
    where the wedge product on the left is the formal wedge product (that is, $\alpha \wedge \beta$ is seen as an element in $\wedge^2V$), while the one on the right is the wedge product of differential forms on $X$, $\alpha \wedge \beta \in V$. 

    This example already leads to an interesting observation: The non-degeneracy of $\omega$ does not impose that $V$ is finite-dimensional in every degree, but only that $V$ is \emph{perfect}. That is, it has finite-dimensional cohomology concentrated in finitely many degrees. \hfill$\triangle$
\end{Example}

\begin{Example}\label{ex2}
    Let $X$ be an $n$-dimensional closed oriented manifold. Let $G$ be a Lie group with Lie algebra $\g$, and $(P, \nabla)$ a principal $G$-bundle with a flat connection $\nabla$. 
    Recall that the adjoint bundle $\ad P := P \times_G\g$ is constructed as the quotient of $P \times \g$ by the action $(p, \xi)\cdot g = (p \cdot g, \Ad_{g^{-1}}(\xi))$. 
    The cochain complex of $\ad P$-valued forms $V = (\Omega^\bullet(X, \ad P), \nabla)[1]$ can be seen as a complex of $\Gamma$-representations, where $\Gamma = \mathrm{Aut}(P, \nabla) \subset C^{\infty}(X, G)$. Recall that a connection $\nabla$ on a vector bundle $E$ defines a map 
    $\Gamma(X,E)=\omega^0(X,E)\to \Omega^1(X,E)$ that can be consistently extended to a degree one endomorphism of $\Omega^\bullet(X,E)$ by the formula $\nabla(\alpha\otimes e)=d_{dR}(\alpha)\otimes e+(-1)^{|\alpha|}\alpha\wedge\nabla(e)$. The connection is flat if and only if this operator squares to zero (i.e.~is a differential). 

    For every $G$-invariant symmetric non-degenerate pairing $\langle -, -\rangle$ on $\g$ we have a $(2-n)$-shifted symplectic structure on $V$:
    \begin{equation*}
        \omega(\alpha \wedge \beta) = \int_X \langle \alpha \wedge \beta \rangle,
    \end{equation*}
    where $\langle \alpha \wedge \beta \rangle$ is the pairing extended to $\ad P$-valued forms. 
\hfill$\triangle$\end{Example}

\subsubsection{The case of $1$-term complexes}

If $V$ is concentrated in degree $d$, then it has the form $V=W[-d]$, for $W$ non-trivial\footnote{The zero cochain complex is $n$-shifted symplectic for every $n$. } and concentrated in degree $0$. Then $V^*$ is concentrated in degree $-d$. By the non-degeneracy condition, $V$ only admits $n$-shifted symplectic structures for $n=-2d$. There are two distinct cases:
\begin{itemize}
    \item If $d$ is odd, then $\wedge^2V = \Sym^2(W)[-2d]$, and an $n$-shifted symplectic structure on $V$ is a scalar product on $W$.
    \item If $d$ is even, then $\wedge^2 V = \wedge^2 W [-2d]$, hence an $n$-shifted symplectic structure on $V$ is an honest symplectic structure on $W$. 
\end{itemize}

\begin{Example}
    Let $G$ be a Lie group. Any non-degenerate $\Ad$-invariant scalar product on $\g$ defines a $2$-shifted symplectic structure on $\g [1]$ as a complex of $G$-representations. 
\hfill$\triangle$\end{Example}

\subsubsection{The case of $2$-term complexes}\label{subsec-2.2.2}

Let $V = (E \overset{a}{\longrightarrow} F)$ be concentrated in degrees $d$ and $d+1$ (more precisely, $V=C[-d-1]$ with $C$ being the cone of $a$). Let's assume that $V$ is not acyclic (that is, $a$ is not an isomorphism), since otherwise $V\simeq 0$ is trivially $n$-shifted symplectic for every $n$ (the $n$-shifted symplectic structure being zero). 

\medskip

If $V$ admits an $n$-shifted symplectic structure and $n$ is even, then either $n=-2d$ and 
\begin{equation*}
   \begin{split}
    &\coker (a) = H^{d+1}(V) \cong H^{-d+1}(V^*)=0,\\
    &\ker (a) = H^d(V) \cong H^{-d}(V^*) = \coker (a^*),\\
    &0=H^{d-1}(V)\cong H^{-d-1}(V^*)=\ker(a^*),
   \end{split}
\end{equation*}
which implies that $V \simeq \ker(a)[-d]$ (meaning in particular that $a$ is surjective), or $n=-2d-2$ and 
\begin{equation*}
   \begin{split}
    &\ker(a) = H^{d}(V) = H^{-d-2}(V^*) = 0,\\
    &\coker (a) = H^{d+1}(V) \cong H^{-d-1}(V^*) = \ker(a^*),\\
    &0 = H^{d+2}(V)\cong H^{-d}(V^*)=\coker(a^*),
   \end{split}
\end{equation*}
which implies that $V \simeq \coker(a)[-d-1]$ (meaning in particular that $a$ is injective). 
In both of these situations, $V$ is quasi-isomorphic to a $1$-term complex. 

Conversely, if $a$ is either surjective or injective, but not an isomorphism, then $n$ must be even. 

\medskip

Let us now assume that $n$ is odd. The condition that $\omega^\flat$ is a quasi-isomorphism imposes that $n=-1-2d$. 
A cochain map $\omega:\wedge^2V\to k[n]$ is completely determined 
by a linear map $\omega_L: E\otimes F \to k$ satisfying the \emph{cochain condition} 
\begin{equation}\label{eq:cochain-condition}
    \omega_L(e_1 \otimes a(e_2)) = (-1)^d \omega_L(e_2\otimes a(e_1)).
\end{equation}
Let us introduce the linear map $\alpha:E\to F^*$ defined by $\alpha(e):=\omega_L(e\otimes-)$. We have the following description of $\omega^\flat$ in terms of $\alpha$: 
\begin{equation}\label{eq-2-term}
V = \left(  
\begin{tikzcd}[row sep=scriptsize]
    E \\
    F
    \arrow["a", from=1-1, to=2-1]
\end{tikzcd} 
\right)
\xrightarrow{\omega^\flat = 
\begin{pmatrix}
    \alpha \\ \alpha^*
\end{pmatrix}}
\left( 
\begin{tikzcd}[row sep=scriptsize]
{F^*} \\
{E^*}
\arrow["{-a^*}", from=1-1, to=2-1]
\end{tikzcd}
\right) 
= V^*[n].
\end{equation}
The cochain condition \eqref{eq:cochain-condition} is equivalent to requiring that $\begin{pmatrix}
    \alpha \\ \alpha^*
\end{pmatrix}$ is a cochain map. The non-degeneracy condition amounts to requiring that $\ker(a) \overset{\alpha}{\longrightarrow} \ker(a^*)$ and $\coker(a) \overset{\alpha^*}{\longrightarrow} \coker(a^*)$  are isomorphisms. One easily sees that $\alpha$ is an isomorphism if and only if $\alpha^*$ is. 

\begin{Remark}
Instead of viewing the commuting square \eqref{eq-2-term} as the cochain map $\omega^\flat$ (going from $V$ to $V^*[n]$), one could view it as the cochain map $\begin{pmatrix}
    a & -a^*
\end{pmatrix}$ going from $(E\overset{\alpha}{\longrightarrow} F^*)$ to $(F\overset{\alpha^*}{\longrightarrow} E^*)$. It turns out that $\omega^\flat$ is a quasi-isomorphism if and only if $\begin{pmatrix}
    a & -a^*
\end{pmatrix}$ is. This is because they both have the same cone. Hence the non-degeneracy condition is equivalent to requiring that $a: \ker(\alpha) \longrightarrow \ker(\alpha^*)$ is an isomorphism. 
\end{Remark}

\begin{Lemma}\label{lem:2-term-cplx-shifted-sympl}
    Assume $E$ and $F$ are finite-dimensional. Then the non-degeneracy condition is equivalent to the condition that $\ker(a)\cap\ker(\alpha)=0$ and $\dim(E) = \dim (F)$. 
\end{Lemma}
\begin{proof}
The proof is a linear algebra exercise which we leave to the reader. 
\end{proof}

\begin{Example}\label{ex-4}
    Let $G$ be a Lie group or an affine algebraic group with Lie algebra $\g$. Let $A = \mathcal{O}(\g^*)$ be the ring of functions on $\g^*$ (in the algebraic case, $\mathcal{O}(\g^*) = \Sym(\g)$). This is a $G$-algebra, because $\g^*$ is a $G$-space with respect to the adjoint action. The action of an arbitrary $g\in G$ on a monomial in $A$ is
    \begin{equation*}
        g \cdot (x^p) := (\Ad_g(x))^p. 
    \end{equation*}
    Consider the infinitesimal action map 
    \begin{equation*}
    \begin{split}
        \g &\longrightarrow \VF(\g^*) \cong A \otimes \g^*\\
        x &\longmapsto \vec{x}.
    \end{split}
    \end{equation*}
    This induces an $A \rtimes G$-module map
    \begin{equation*}
    \begin{split}
        a: A \otimes \g &\longrightarrow \VF(\g^*)\cong A \otimes \g^*\\
        f \otimes x &\longmapsto f\vec{x}.
    \end{split}
    \end{equation*}
    We view this as a $2$-term complex of $A\rtimes G$-modules concentrated in degrees $-1$ and $0$. Then we define 
    \begin{equation}
        \omega_L: (A\otimes \g) \otimes_A (A\otimes \g^*) \cong A \otimes \g \otimes \g^* \xrightarrow{id \otimes ev} A.
    \end{equation}
    Let us check the cochain condition\eqref{eq:cochain-condition} with $d=-1$: if  $(e_i)_{i=1,\dots, n}$ is a basis of $\g$, then $a(e_i) = c_{ij}^k e_j^*e_k$, and thus 
    \begin{equation*}
        ev(e_i,a(e_j))=ev(e_i,c_{js}^ke_s^*e_k)=c_{ji}^ke_k=-c_{ij}^ke_k=-ev(e_j,a(e_i)).
    \end{equation*}
    This induces an isomorphism of $2$-term complexes
    \begin{equation*}
    \left(  
    \begin{tikzcd}[row sep=scriptsize]
        {A \otimes \g} \\
        {A \otimes \g^*}
        \arrow["a", from=1-1, to=2-1]
    \end{tikzcd} 
    \right)
    \overset{id}{\longrightarrow}
    \left( 
    \begin{tikzcd}[row sep=scriptsize]
	{A \otimes \g} \\
	{A \otimes \g^*}
	\arrow["{-a^*}", from=1-1, to=2-1]
    \end{tikzcd}
    \right).
    \end{equation*}
    One can indeed check that $a=-a^*$ on basis elements:
    \begin{equation*}
        a(e_i)(e_j)
        =  c_{ij}^k e_k
        = -c_{ji}^k e_k
        = -a(e_j)(e_i)
        = - a^*(e_i)(e_j).
    \end{equation*}
    Therefore, this defines a $1$-shifted symplectic structure on the $2$-term complex $A\otimes \g \xrightarrow{a} \VF(\g^*)\cong A\otimes \g^*$.
\hfill$\triangle$\end{Example}

\begin{Example}
    Let $G$ be a reductive algebraic group over $k$ (or a compact group when $k=\mathbb{R}$) with a choice of a non-degenerate invariant pairing $\langle \cdot, \cdot\rangle:\Sym^2(\g)\to k$ on $\g$. The ring of functions  
    $B = \mathcal{O}(G)$ is a $G$-algebra because $G$ is a $G$-space with respect to the conjugation action. 
    Consider the infinitesimal action map 
    \begin{equation*}
    \begin{split}
        a: \g &\longrightarrow \VF(G)\\
        x &\longmapsto \vec{x} = x^L-x^R.
    \end{split}
    \end{equation*}
    By choosing the left trivialization $TG \cong G \times \g$ we get an isomorphism $\VF(G)\cong B\otimes \g$, and $a(x)_g=x-\Ad_g(x)$ for every $x\in\g$ and every $g\in G$. 
    This gives us a $2$-term complex of $B\rtimes G$-modules 
    \begin{equation*}
    \begin{split}
        B\otimes \g &\overset{a}{\longrightarrow} \VF(G)\\
        f\otimes x &\longmapsto f\vec{x},
    \end{split}
    \end{equation*}
    concentrated in degrees $-1$ and 0 (left to right).
    Let us now define 
    \begin{equation*}
        \omega_L: (B\otimes \g)\otimes_B\VF(G) \cong (B\otimes \g)\otimes_B (B\otimes \g)\cong B\otimes \g\otimes \g \xrightarrow{id\otimes \langle \cdot, \cdot\rangle} B,
    \end{equation*}
    where the first map is the average of left and right Maurer-Cartan forms
    \begin{equation*}
        \frac{1}{2}(g^{-1}dg + dg g^{-1}) \in (\g \otimes \Omega^1(G))^G.
    \end{equation*}

    \begin{Exercise} Check that the induced map $\alpha: B\otimes \g \to \Omega^1(G) \cong \VF(G)$ is given by $\alpha(1 \otimes x) = \frac{1}{2}(x^L + x^R)$, after identifying $1$-forms and vector fields through the pairing: $\Omega^1(G) \cong B\otimes \g^* \cong B \otimes \g \cong \VF(G)$. 
\end{Exercise}
    Observe now that $\rk_B(B\otimes \g) = \dim (\g ) = \rk_B(\Omega^1(G))$ and that 
    \begin{equation*}
    \ker(a)\cap\ker(\alpha) = \{x\mid x^L-x^R = 0 = \frac{1}{2}(x^L + x^R)\} = 0.
    \end{equation*}
    Then, by \cref{lem:2-term-cplx-shifted-sympl}, $\omega_L$ is a $1$-shifted symplectic structure. 
\hfill$\triangle$\end{Example}

\begin{Example}
    The previous examples are specific cases of the following. Let $G_\bullet$ be a Lie groupoid. (See for example \cite{CdSWeinstein, Mackenzie, MoerdijkMrcun}).
    This can be represented by a diagram 
\[\begin{tikzcd}
	{G_1} & {G_0}
	\arrow["s"{description}, shift left=2, from=1-1, to=1-2]
	\arrow["t"{description}, shift right, from=1-1, to=1-2]
	\arrow["e"{description}, curve={height=-12pt}, from=1-2, to=1-1]
\end{tikzcd}\]
    where $G_0$ is the manifold of objects, $G_1$ is the manifold of arrows, and the two form a category where 
    \begin{itemize}
        \item The source and target maps $s,t$ are surjective submersions.
        \item The map $e$ is the map associating to each object its identity arrow. 
        \item All the arrows in $G_1$ are invertible. 
    \end{itemize}  
    The space of composable arrows of $G_\bullet$ is
    \begin{equation*}
        G_2 = G_1 \underset{t,G_0,s}{\times} G_1 
        = \{(h,g)\in (G_1)^2 \mid sh = tg\},
    \end{equation*}
    and the groupoid multiplication (i.e. the composition, when considered as a category), is the map 
    \begin{equation*}
    \begin{split}
    m: G_2 &\longrightarrow G_1\\
    (h,g) &\longmapsto hg.
    \end{split}
    \end{equation*}

    Consider a \emph{multiplicative} $2$-form $\omega \in \Omega^2(G_1)$. This is a $2$-form such that 
    \begin{equation}\label{eq:omega-multiplicative}
        pr_1^*\omega - m^*\omega + pr_2^*\omega = 0 \in \Omega^2(G^{(2)}).
    \end{equation}

    Let $L$ be the Lie algebroid associated to $G$. (See for example \cite{CdSWeinstein, Mackenzie, MoerdijkMrcun}).
    As a vector bundle, this is
    \begin{equation*}
        L := e^*T^sG_1 = e^*\ker(Ts: TG_1 \to s^*TG_0) 
        \longrightarrow G_0, 
    \end{equation*}
    the tangent space to the source fibers at the unit section.
    The anchor map is the differential of the target:
    \begin{equation*}
        a = Tt: L \to TG_0. 
    \end{equation*}
    This can be seen as a $G_\bullet$-equivariant $2$-term complex of bundles over $G_0$ (concentrated in degrees $-1$ and $0$). 

    We define $\omega_L$ to be the restriction of $\omega$ to 
    \begin{equation*}
        L \times TG_0 \subset e^*(TG_1 \times TG_1).
    \end{equation*}
    One can check that the multiplicativity of $\omega$ as in \eqref{eq:omega-multiplicative} implies the fact that $\omega_L$ satisfies the cochain condition \eqref{eq:cochain-condition}. 

    In \cite{Xu}, $\omega$ is called an \emph{almost quasi-symplectic structure} on $G_\bullet$ if 
    \begin{equation*}
        a: \ker(\alpha_x) \to \ker(\alpha_x^*)
    \end{equation*}
    is an isomorphism for every $x\in G_0$,
    i.e. if $\omega_L$ defines a $1$-shifted symplectic structure on the $2$-term complex $L\overset{a}{\longrightarrow} TG_0$ of 
    ($G_\bullet$-equivariant) vector bundles. In other words, for every $x\in G_0$,  $\omega_{L,x}$ defines a $1$-shifted symplectic 
    structure on the $2$-term complex $L_x\overset{a_x}{\longrightarrow} T_xG_0$. 
\hfill$\triangle$\end{Example}

\begin{Example}\label{ex:-1-shifted-structure-lagrangian-composition}
    Let $(V, \omega)$ be a usual symplectic vector space and $L_1,L_2\subset V$ lagrangian subspaces in the usual sense. 
    Consider the $2$-term complex 
    \begin{equation*}
    \begin{split}
        L_1 \oplus L_2 &\longrightarrow V\\
        (\ell_1,\ell_2) &\longmapsto \ell_1-\ell_2,
    \end{split}
    \end{equation*}
    with $L_1 \oplus L_2$ in degree 0 and $V$ in degree 1. 
    Define 
    \begin{equation*}
    \begin{split}
        \omega_L:(L_1 \oplus L_2) \otimes V &\longrightarrow k\\
        (\ell_1,\ell_2) \otimes v &\longrightarrow \omega(\ell_1+\ell_2,v).
    \end{split}
    \end{equation*}
    This satisfies the cochain condition \eqref{eq:cochain-condition} with $d=0$:
    \begin{equation*}
        \omega(\ell_1+\ell_2, \ell_1'-\ell_2') 
        -\omega(\ell_1'+ \ell_2', \ell_1 - \ell_2) 
        =2\omega(\ell_1,\ell_1') - 2\omega(\ell_2,\ell_2') 
        = 0,
    \end{equation*}
    because $L_1$ and $L_2$ are isotropic. 
    It also satisfies the non-degeneracy condition, since $\dim (L_1\oplus L_2) = \dim (V)$ and 
    \begin{align*}
        \ker(a)\cap\ker(\alpha)&=\{(\ell_1,\ell_2)\in L_1\oplus L_2 \mid \ell_1=\ell_2, \text{ and } \omega(\ell_1+\ell_2, -) = 0 \}\\& \cong \{\ell\in L_1\cap L_2\mid\omega(\ell,-)=0\} = 0,
    \end{align*}
    where the last equality holds because $\omega$ is non-degenerate. 
    Therefore, this defines a $(-1)$-shifted symplectic structure. Notably, observe that we used all of the assumptions to show this. 
\hfill$\triangle$\end{Example}

\subsection{Lagrangian structures}

\subsubsection{Recollection about homotopies}

We begin by recalling the notion of homotopy and introducing the notion of \emph{cocone} of a cochain map. This is also called \emph{mapping cocone}, or \emph{homotopy fiber}.

\begin{Definition}
    Consider two cochain maps
    \begin{equation*}
        \phi, \psi: (V,\delta_V) \longrightarrow (W,\delta_W).
    \end{equation*}
    A \emph{homotopy} $\eta$ between $\phi$ and $\psi$ is a map of graded vector spaces
    \begin{equation*}
        \eta: V \longrightarrow W[-1]
    \end{equation*}
    such that
    \begin{equation*}
        \eta \delta_V + \delta_W \eta = \phi - \psi.
    \end{equation*}
    In this case we write $\phi \overset{\eta}{\sim} \psi$.
\end{Definition}

There is a natural composition for homotopies, given by the sum: 
if $\phi \overset{\eta}{\sim} \psi$ and $\psi \overset{\nu}{\sim} \kappa$ then $\phi\overset{\eta+\nu}{\sim} \kappa$. 

\begin{Definition}
    The \emph{cocone} of a cochain map
    \begin{equation*}
        \phi:(V, \delta_V) \longrightarrow (W, \delta_W)
    \end{equation*}
    is the cochain complex
    \begin{equation*}
        \hofib(\phi)=\left(\begin{array}{c}
            V\\
            \oplus\\
            W[-1]
        \end{array}
        , \delta
        \right),
    \end{equation*}
    with differential 
    \begin{equation*}
        \delta = \begin{pmatrix}
            \delta_V &0 \\
            \phi &-\delta_W
        \end{pmatrix}.
    \end{equation*} 
\end{Definition}

The main property of this construction is the fact that a cochain map 
\begin{equation*}
    (U, \delta_U) \longrightarrow \hofib(\phi)
\end{equation*}
coincides with the data of a cochain map
\begin{equation*}
    \psi: (U, \delta_U) \longrightarrow (V,\delta_V),
\end{equation*}
together with a homotopy 
\begin{equation*}
    \phi\psi \overset{\eta}{\sim} 0,
\end{equation*}
and we write
\begin{equation*}
    \begin{pmatrix}
        \psi \\
        \eta 
    \end{pmatrix}: (U, \delta_U) \longrightarrow \hofib(\phi).
\end{equation*}

Indeed, after denoting the cocone by $(C, \delta_C):=\hofib(\phi)$, we have that 
\begin{equation*}
\begin{split}
    \begin{pmatrix}
        \psi \delta_C\\
        \eta \delta_C
    \end{pmatrix}
    = 
    \begin{pmatrix}
        \psi \\
        \eta 
    \end{pmatrix}
    \delta_C
    &= 
    \delta
    \begin{pmatrix}
        \psi \\
        \eta 
    \end{pmatrix}\\
    &=
    \begin{pmatrix}
        \delta_V  &0\\
        \phi &-\delta_W
    \end{pmatrix}
    \begin{pmatrix}
        \psi \\
        \eta 
    \end{pmatrix}
    = 
    \begin{pmatrix}
        \delta_v \psi \\
        \phi\psi - \delta_W \eta
    \end{pmatrix}.
\end{split}
\end{equation*}

\subsubsection{Isotropic structures}

\begin{Definition}
    Let $(V, \omega)$ be a complex together with a $2$-form $\omega:\wedge^2 V \to k[n]$.
    Let $\phi:L \to V$ be a cochain map. An \emph{isotropic structure} on $\phi$ (w.r.t. $\omega$) is a homotopy $\omega|_L \overset{\eta}{\sim} 0$. Here $\omega|_L$ is an abbreviated notation for $\omega(\wedge^2\phi)$. 
   \end{Definition}
 Concretely, this is a map
    \begin{equation*}
        \eta: \wedge^2 L \to k[n-1],
    \end{equation*}
    such that 
    \begin{equation*}
        \eta(\delta a\wedge b) 
        +(-1)^{|a|} 
        \eta(a \wedge \delta b) 
        = \omega(\phi(a) \wedge \phi(b)),
    \end{equation*}
    for any $a,b \in L$. 

\begin{Remark}
    Note that the above definition implies that, if $\eta$ is an isotropic structure on $\phi$ with respect to $\omega$, then $\eta^\flat$ provides a homotopy between $\phi^*\omega^\flat \phi:L\to L^*[n]$ and $0$ (as such, $\eta^\flat$ is in particular a graded map $L\to L^*[n-1]$). 
    Therefore, we have a morphism of complexes 
    \begin{equation*}
    \begin{pmatrix}
        \phi\\
        \eta^\flat
    \end{pmatrix}
    : L \longrightarrow 
    \hofib\big(\phi^*\omega^\flat:V\to L^*[n]\big).
    \end{equation*}
\end{Remark}

\subsubsection{Non-degeneracy condition}

\begin{Definition}\label{def-nd-homotopy}
    Borrowing the notation from the previous subsection, we say that $\eta$ is \emph{non-degenerate} if $\begin{pmatrix}
        \phi\\
        \eta^\flat
    \end{pmatrix}$ is a quasi-isomorphism.
 A non-degenerate isotropic structure is called a \emph{lagrangian structure}.
\end{Definition}

\begin{Remark}
    Alternatively, we could have asked that 
    \begin{equation*}
        \begin{pmatrix}
            \omega^\flat\phi\\
            \eta^\flat
        \end{pmatrix}:
        L \to \hofib(\phi^*)
    \end{equation*}
    be a quasi-isomorphism, as this is actually equivalent.  
    As a consequence, we have a morphism between long exact sequences 
\[\begin{tikzcd}
	\dots & {H^\bullet(L)} & {H^\bullet(V)} & {H^\bullet(L^*[n])} & {H^{\bullet + 1}(L)} & \dots \\
	\dots & {H^\bullet(L)} & {H^\bullet(V^*[n])} & {H^\bullet(L^*[n])} & {H^{\bullet + 1}(L)} & \dots
	\arrow[from=1-1, to=1-2]
	\arrow[from=1-2, to=1-3]
	\arrow[equals, from=1-2, to=2-2]
	\arrow[from=1-3, to=1-4]
	\arrow[from=1-3, to=2-3]
	\arrow[from=1-4, to=1-5]
	\arrow[equals, from=1-4, to=2-4]
	\arrow[from=1-5, to=1-6]
	\arrow[equals, from=1-5, to=2-5]
	\arrow[from=2-1, to=2-2]
	\arrow[from=2-2, to=2-3]
	\arrow[from=2-3, to=2-4]
	\arrow[from=2-4, to=2-5]
	\arrow[from=2-5, to=2-6]
\end{tikzcd}\]
Therefore, we get that 
\begin{equation*}
    H^\bullet (V) 
    \longrightarrow H^\bullet (V^*[n]),
\end{equation*}
is an isomorphism, i.e. $\omega$ is non-degenerate.
\end{Remark}

\begin{Example}
    Let $Y$ be an $(n+1)$-dimensional compact oriented manifold with boundary $\partial Y = X$. 
    Consider the cochain complex $V=(\Omega^\bullet(X), d_{dR})[1]$ equipped with the $(2-n)$-shifted symplectic structure 
    \begin{equation*}
        \omega (\alpha \wedge \beta) := \int_X \alpha \wedge \beta,
    \end{equation*}
    from \cref{ex:shift-sympl-struct-deRham-cplx}.
    The restriction of a form on $Y$ to a form on the boundary $X$ gives a cochain map
    \begin{equation*}
    \begin{split}
        \phi: 
        L = (\Omega^{\bullet}(Y), d_{dR}) 
        &\longrightarrow
        (\Omega^{\bullet}(X), d_{dR}),\\
        \alpha
        &\longrightarrow
        \alpha|_X.
    \end{split}
    \end{equation*}
    We claim that the map $\eta: \wedge^2 L \to \R [n-1]$, defined by 
    \begin{equation*}
        \eta(\alpha \wedge \beta) 
        := \int_Y \alpha \wedge \beta,
    \end{equation*}
    defines a lagrangian structure on $\phi$ with respect to $\omega$. 
    Firstly, it is isotropic, i.e. a homotopy $\omega|_L \sim 0$, by Stokes' theorem:\footnote{We do the following computation within $\Omega^\bullet(Y)$, i.e.~without the degree shift, in order to ease the presentation and remove a few minus signs here and there. This does not affect the result, as the sign change due to the degree shift $|\alpha|_L=|\alpha|-1$ is compensated by the sign modification $d_L=-d_{dR}$ for the differential. }
    \begin{equation*}
    \begin{split}
    \eta(d_{dR}\alpha \wedge \beta +(-1)^{|\alpha|} \alpha \wedge d_{dR}\beta) 
    &= \int_Y d_{dR}\alpha \wedge \beta +(-1)^{|\alpha|} \alpha \wedge d_{dR}\beta\\
    &= \int_Y d_{dR}(\alpha \wedge \beta) 
    = \int_X (\alpha \wedge \beta)|_X 
    = \omega(\alpha|_X \wedge \beta|_X).
    \end{split}
    \end{equation*}
    Secondly, it is non-degenerate:
    Consider $\phi^*:V^*[n] \to L^*[n]$. The cohomology of $\hofib(\phi^*)$ is dual (up to a shift) to the relative cohomology $H^\bullet(Y,X)$, that is itself dual (up to the same shift) to $H^\bullet(Y) = H^{\bullet - 1}(L)$. 
    Non-degeneracy thus follows from relative Poincaré duality. 
\hfill$\triangle$\end{Example}

\begin{Example}
    The same example can be repeated with the following additional structure, as in \cref{ex2}: let $\g$ be a Lie algebra with an invariant non-degenerate pairing and $(P,\nabla)$ a flat principal $G$-bundle. Then consider the following cochain complexes:
    \begin{equation*}
    V=(\Omega^\bullet(X, \ad(P|_X)), \nabla)[1], 
    \quad
    L=(\Omega^\bullet(Y, \ad(P)), \nabla)[1]. 
    \end{equation*}
\hfill$\triangle$\end{Example}

Recall Weinstein's ``symplectic creed'' \cite{WeinsteinCreed}: \emph{everything is a lagrangian submanifold}. As one can see from the examples above, not all lagrangian morphisms are subcomplexes, hence one could be tempted to modify the symplectic creed as follows: \emph{everything is a lagrangian morphism/structure}. 

\begin{Example}
A surprising incarnation of this more ``homotopical'' symplectic creed is that a symplectic structure is a particular example of a lagrangian structure. Let $V=0$, equipped with the zero $n$-shifted symplectic structure $\omega=0$. Then a lagrangian structure on $L \to 0$ is a map
\begin{equation*}
    \eta: \wedge^2 L \to k[n-1]
\end{equation*}
such that 
\begin{equation*}
    \eta(\delta a \wedge b) +(-1)^{|a|} \eta(a \wedge \delta b) = 0|_L=0, 
\end{equation*}
that is non-degenerate. 
In particular, we see that $\eta$ is a degree $n-1$ skew-symmetric pairing on $L$, while the non-degeneracy condition implies that the cochain map 
\begin{equation*}
    \eta^\flat=\begin{pmatrix}
        0\\
        \eta^\flat
    \end{pmatrix}
    : L \longrightarrow \hofib(0 \to L^*[n]) = L^*[n-1]
\end{equation*}
is a quasi-isomorphism; i.e. $\eta$ defines an $(n-1)$-shifted symplectic structure on $L$. Running the same reasoning backwards, one conversely gets that any $(n-1)$-shifted symplectic structure $\omega$ on $V$ defines a lagrangian structure on $V\to 0$ (which is $\omega$ itself, now viewed as a self-homotopy of zero).  
\hfill$\triangle$\end{Example}

\subsubsection{Back to $2$-term complexes}\label{subsec-2.3.4}

Let us consider again the case of a two-term complex 
\begin{equation*}
    V=\left( E \overset{a}{\longrightarrow} F \right),
\end{equation*}
with $E$ in degree $d$ and $F$ in degree $d+1$, equipped with an $n$-shifted symplectic structure $\omega$ for $n=-1-2d$. Recall that $\omega$ is determined by $\omega_L:E \otimes F \to k$, and that we write $\alpha:=\omega_L^\flat: E \to F^*$. 

Let now $L$ be the two-term complex
\begin{equation*}
    L = \left( E \overset{b}{\longrightarrow} B \right),
\end{equation*}
with $E$ in degree $d$ and $B$ in degree $d+1$. Consider a cochain map 
\begin{equation*}
    \phi: L \longrightarrow V.
\end{equation*}
given as the matrix 
\begin{equation*}
    \phi = \begin{pmatrix}
        id_E\\
        f
    \end{pmatrix},
\end{equation*}
for some map $f:B \longrightarrow F$. In order for $\phi$ to be a cochain map one must require that $fb=a$. 

We would like to answer the following question:
\begin{center}
    \emph{What is a lagrangian structure on $\phi$ with respect to $\omega$?}
\end{center}

First of all, an isotropic structure on $\phi$ is a graded map $\eta: \wedge^2 B \to k$ (for degree reasons) such that for any $x \in E$ and $y\in B$, 
\begin{equation*}
    \eta(bx\wedge y) = \omega_L(x\wedge f(y))
\end{equation*}
(this is the property of being a homotopy between $\omega$ and $0$). 
This last condition is equivalent to 
\begin{equation*}
    \eta^\flat \circ b = f^* \circ \alpha,
\end{equation*}
and thus to 
\begin{equation*}
    b^*\circ \eta^\flat = \alpha^* \circ f.
\end{equation*}
In other words, $\eta^\flat$ provides a homotopy between the composed map 
\[\begin{tikzcd}
	L & V & {V^*[n]} & {L^*[n],}
	\arrow["\phi", from=1-1, to=1-2]
	\arrow["{\omega^\flat}", from=1-2, to=1-3]
	\arrow["{\phi^*}", from=1-3, to=1-4]
\end{tikzcd}\]
which can be written more precisely as the composition of the horizontal maps 
\[\begin{tikzcd}
	E & E & {F^*} & {B^*} \\
	B & F & {E^*} & {E^*}
	\arrow["{id_E}", from=1-1, to=1-2]
	\arrow["b", from=1-1, to=2-1]
	\arrow["\alpha", from=1-2, to=1-3]
	\arrow["a", from=1-2, to=2-2]
	\arrow["{f^*}", from=1-3, to=1-4]
	\arrow["{a^*}", from=1-3, to=2-3]
	\arrow["{b^*}", from=1-4, to=2-4]
	\arrow["f"', from=2-1, to=2-2]
	\arrow["{\alpha^*}"', from=2-2, to=2-3]
	\arrow["{id_{E^*}}"', from=2-3, to=2-4]
\end{tikzcd}\]
between two-term complexes, and the zero map $L \overset{0}{\longrightarrow} L^*[n]$.

Then, it remains to express the non-degeneracy condition for such an isotropic structure on $\phi$ in terms of $\eta$. The cochain map 
\begin{equation*}
    \begin{pmatrix}
        \phi\\
        \eta^\flat
    \end{pmatrix}:
    L \longrightarrow \hofib(\phi^* \circ \omega^\flat)
\end{equation*}
is the horizontal map in the diagram
\begin{equation}\label{eq-3-term}
    \begin{tikzcd}
	{\text{degree}\quad  d} & E & E \\
	{\text{degree}\quad  d+1} & B & {F\oplus B^*} \\
	{\text{degree}\quad  d+2} & 0 & {E^*.}
	\arrow["{id_E}", from=1-2, to=1-3]
	\arrow["b", from=1-2, to=2-2]
	\arrow["{a \oplus (f^*\circ \alpha)}", from=1-3, to=2-3]
	\arrow["{f\oplus\eta^\flat}"', from=2-2, to=2-3]
	\arrow[from=2-2, to=3-2]
	\arrow["{(\alpha^*, -b^*)}", from=2-3, to=3-3]
	\arrow["0"', from=3-2, to=3-3]
\end{tikzcd}
\end{equation}
Recall that $\alpha^*: F \to E^*$ is an isomorphism (because $\omega$ is non-degenerate by definition). Therefore the projection 
\[\begin{tikzcd}
	E & E \\
	{F\oplus B^*} & {B^*} \\
	{E^*} & 0
	\arrow["{id_E}", from=1-1, to=1-2]
	\arrow["{a \oplus (f^*\circ \alpha)}"', from=1-1, to=2-1]
	\arrow["{f^*\circ\alpha}",from=1-2, to=2-2]
	\arrow[two heads, from=2-1, to=2-2]
	\arrow["{(\alpha^*, -b^*)}"', from=2-1, to=3-1]
	\arrow[from=2-2, to=3-2]
	\arrow[two heads, from=3-1, to=3-2]
\end{tikzcd}\]
is a quasi-isomorphism. 
As a result, \eqref{eq-3-term} is a quasi-isomorphism (meaning that the isotropic structure is non-degenerate) if and only if 
\[\begin{tikzcd}
	E & E \\
	B & B^*
	\arrow["{id_E}", from=1-1, to=1-2]
	\arrow["b", from=1-1, to=2-1]
	\arrow["{f^*\circ \alpha}", from=1-2, to=2-2]
	\arrow["{\eta^\flat}"', from=2-1, to=2-2]
\end{tikzcd}\]
is a quasi-isomorphism, which in turn is true if and only if $\eta^\flat$ is an isomorphism.

\begin{Example}
    Recall the $1$-shifted symplectic structure on the $2$-term complex 
    \begin{equation*}
        \mathcal{O}(\g^*)\otimes \g \longrightarrow\VF(\g^*)
    \end{equation*}
    from \cref{ex-4}. 
    Let $X$ be a smooth affine algebraic variety or a smooth manifold, and let $\mu:X \to \g^*$ be a $G$-equivariant map. Then the functor 
    \begin{equation*}
        \mu^*=\mathcal O(X)\otimes_{\mathcal O(\g^*)}-: \mathcal{O}(\g^*)\text{-}\mathsf{Mod} \longrightarrow \mathcal{O}(X)\text{-}\mathsf{Mod}
    \end{equation*}
    is symmetric monoidal\footnote{Thus the image by $\mu^*$ of an object equipped with an $n$-shifted skew-symmetric pairing is also equipped with an $n$-shifted skew-symmetric pairing. } and exact\footnote{Hence the non-degeneracy property  is preserved under $\mu^*$. }. 
    Therefore the complex 
    \begin{equation*}
        V:= \left(
        \mathcal{O}(X)\otimes \g
        \overset{a}{\longrightarrow} 
        \Gamma(X, \mu^*T\g)
        \right)
    \end{equation*}
    carries a $1$-shifted symplectic structure (in $\mathcal O(X)\rtimes G$-modules). 
    Let $L$ be the complex
    \begin{equation*}
        L:= \left(
        \mathcal{O}(X)\otimes \g
        \overset{b}{\longrightarrow} 
        \VF(X)
        \right),
    \end{equation*}
    where $b$ is the infinitesimal action, and consider the map $L\overset{\phi}{\longrightarrow} V$ given by the identity on $\mathcal{O}(X)\otimes \g$ and $f=\mu_*: \VF(X) \longrightarrow \Gamma(X, \mu^*T\g)$.

    An isotropic structure on $\phi$ with respect to $\omega$ is a morphism 
    \begin{equation*}
       \eta: \wedge^2_{\mathcal{O}(X)} \VF(X) \longrightarrow \mathcal{O}(X)
    \end{equation*}
    of $\mathcal{O}(X) \rtimes G$-modules, that is to say a $G$-invariant $2$-form $\eta\in \Omega^2(X)^G$, such that, for any $v \in \g$ seen as an element of $\mathcal{O}(\g^*) \otimes \g \cong \Omega^1(\g^*)$, 
    \begin{equation*}
        \iota_{\Vec{v}_X}\eta=\eta(b(v) \wedge - ) = f^*\alpha(v)=\mu^*v. 
    \end{equation*}
    This coincides with requiring that $\mu$ is a moment map for $\eta$. 
    Such an isotropic structure is non-degenerate if and only if $\eta^\flat$ is an isomorphism, since $\alpha$ is an isomorphism in this case. That is, when $\eta$ is almost symplectic. 
\hfill$\triangle$\end{Example}

\begin{Remark}
With a bit more work, one can prove that generalized moment maps, such as the Lie group valued moment maps of Alekseev--Malkin--Meinrenken \cite{AMM} and more generally, moment maps with values in (almost) quasi-symplectic groupoids \cite{Xu}, define Lagrangian structures. 
\end{Remark}

\subsection{Lagrangian correspondences}

\begin{Definition}
    A \emph{Lagrangian correspondence} is the data of two $n$-shifted symplectic complexes $(V_1, \omega_1)$ and $(V_2, \omega_2)$, together with a lagrangian morphism 
    \begin{equation*}
        L \to (V_1 \oplus V_2, \omega_1 - \omega_2). 
    \end{equation*}
\end{Definition}

Lagrangian correspondences have a well-defined composition, provided that they are composed by using \emph{homotopy fiber products}. 

\begin{Definition}
    Let $L_{12}$ be a lagrangian correspondence between $V_1$ and $V_2$ and $L_{23}$ be a lagrangian correspondence between $V_2$ and $V_3$. Then we define their \emph{composition} $L_{13}$ as a lagrangian correspondence between $V_1$ and $V_3$ by taking 
    \begin{equation*}
        L_{13}:= \hofib(L_{12}\oplus L_{23} \xrightarrow{a_2- b_2} V_2), 
    \end{equation*}
    that is $L_{12}\oplus L_{23} \oplus V_2[-1]$ with differential 
    \begin{equation*}
        \begin{pmatrix}
            \delta_{L_{12}} &0 &0\\
            0 &\delta_{L_{23}} &0\\
            a_2 &-b_2 &-\delta_{V_2}
        \end{pmatrix},
    \end{equation*}
    for $a_2, b_2$ as in the diagram
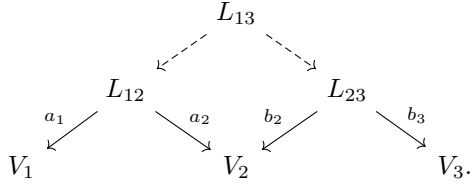
\[\begin{tikzcd}[sep=scriptsize]
	&& {L_{13}} \\
	& {L_{12}} && {L_{23}} \\
	{V_1} && {V_2} && {V_3.}
	\arrow["{p}"', dashed, from=1-3, to=2-2]
	\arrow["{q}", dashed, from=1-3, to=2-4]
	\arrow["{a_1}"', from=2-2, to=3-1]
	\arrow["{a_2}", from=2-2, to=3-3]
	\arrow["{b_2}"', from=2-4, to=3-3]
	\arrow["{b_3}", from=2-4, to=3-5]
\end{tikzcd}\]
\end{Definition}

The lagrangian structure on $L_{13} \to (V_1 \oplus V_3, \omega_1 - \omega_3)$ is defined in the following way.

Recall that the lagrangian structure on $L_{12} \to (V_1 \oplus V_2, \omega_1-\omega_2)$ is a homotopy
\begin{equation*}
    \eta_{12}:\wedge^2 L_{12} \longrightarrow k[n-1]
\end{equation*}
between $a_1^*\omega_1$ and $a_2^*\omega_2$. We also have a homotopy
\begin{equation*}
    \eta_{23}:\wedge^2 L_{23} \longrightarrow k[n-1]
\end{equation*}
between $b_2^*\omega_2$ and $b_3^*\omega_3$.
Let $p:L_{13} \to L_{12}$ and $q:L_{13} \to L_{23}$ be the two canonical projections. 
Observe that $a_2\circ p$ and $b_2 \circ q$ are homotopic: the homotopy is given by by the projection 
\begin{equation*}
    h: L_{13} \longrightarrow V_2[-1],
\end{equation*}
and this induces a homotopy 
\[
h_{\omega_2}:=\frac12\omega_2\big((a_2p+b_2q)\otimes h+h\otimes(a_2p+b_2q)\big)
\]
between $p^*a_2^*\omega_2$ and $q^*b_2^*\omega_2$, due to the following: 
\begin{Lemma}
    Let $L,W$ be complexes, $\omega:\wedge^2W\to k[n]$ and $f,g:L\to W$ be cochain maps, and $h:L\to W[-1]$ be a homotopy between $f$ and $g$. Then 
    \[
h_{\omega}:=\frac12\omega\big(h\otimes(f+g)+(f+g)\otimes h\big)
\]
is a homotopy between $f^*\omega$ and $g^*\omega$. 
\end{Lemma}
\begin{proof}
First observe that 
\[
f^*\omega-g^*\omega=\omega(f\otimes f-g\otimes g)
=\frac12\omega\big((f-g)\otimes(f+g)+(f+g)\otimes(f-g)\big).
\]
Therefore, 
\begin{align*}
h\delta_{\wedge^2L}+\delta_{k[n]}h
& =h_\omega(\delta_L\otimes id_L+id_L\otimes \delta_L) \\ 
& = \frac12\omega\big(h\delta\otimes(f+g)-(f+g)\delta\otimes h+h\otimes(f+g)\delta+(f+g)\otimes h\delta\big) \\
& = \frac12\omega\big(h\delta\otimes(f+g)-\delta(f+g)\otimes h+h\otimes \delta(f+g)+(f+g)\otimes h\delta\big) \\
& = \frac12\omega\big(h\delta\otimes(f+g)+(f+g)\otimes \delta h+\delta h\otimes (f+g)+(f+g)\otimes h\delta\big) \\& = \frac12\omega\big(\underbrace{(h\delta+\delta h)}_{f-g}\otimes(f+g)+(f+g)\otimes \underbrace{(h\delta+\delta h)}_{f-g}\big)
\end{align*}
coincides with $f^*\omega-g^*\omega$. 
\end{proof}
Consequently, $
    p^*\eta_{12} + h_{\omega_2} + q^*\eta_{23}$
defines a homotopy between $p^*a_1^*\omega_1$ and $q^*b_3^*\omega_3$.

\begin{Exercise}
Prove that this isotropic structure is non-degenerate whenever $\eta_{12}$ and $\eta_{23}$ are. 
\end{Exercise}

\begin{Example}[$0$-shifted lagrangian intersection]
    Let $V_1 =V_3=0$, and let $V_2=V$ be an honest symplectic vector space (i.e.~$0$-shifted symplectic). Let $L_{12}$ and $L_{23}$ be honest lagrangian subspaces in $V_2$. Then their composition is the complex
    \begin{equation*}
    \begin{split}
        L_{13} = 
        (L_{12} \oplus L_{23} &\longrightarrow V_2)\\
        (a,b) &\longmapsto a-b, 
    \end{split}
    \end{equation*}
    and $L_{13}\to 0$ carries a Lagrangian structure. Therefore, $L_{13}$ is $(-1)$-shifted symplectic. In the above notation $\eta_{12} =0$ and $\eta_{23} =0$, so 
    $
    \eta = h_\omega$
    is a $(-1)$-shifted symplectic structure, which coincides with (a half of) the one from \cref{ex:-1-shifted-structure-lagrangian-composition}: 
    one can easily check that $\eta_L((a,b),v)=\frac12\omega(a+b,v)$. 
\hfill$\triangle$\end{Example}

\begin{Example}[composition of transverse $0$-shifted lagrangian correspondences]
Let $V_1,V_2,V_3$ be genuine symplectic vector spaces, and assume that $L_{12}$ and $L_{23}$ are genuine lagrangian subspaces (in 
$(V_1 \oplus V_2, \omega_1-\omega_2)$ and $(V_2 \oplus V_3, \omega_2-\omega_3)$, respectively). If $L_{12}$ and $L_{23}$ intersect transversally, 
meaning that  the map 
\[
L_{12}\oplus L_{23}\longrightarrow V_2
\]
is surjective, then $L_{13}=\hofib(L_{12}\oplus L_{23}\to V_2)\simeq\ker(L_{12}\oplus L_{23}\to V_2)$ is thus the usual (strict) fiber product. 
We recover that the usual fiber product of two transversally intersecting usual lagrangian correspondences is again lagrangian. 

Actually, if the intersection is not transversal, then 
$L_{13}\simeq U\oplus W[-1]$, with $U=\ker(L_{12}\oplus L_{23}\to V_2)$ the usual fiber product and $W\neq0$. Observe that $U\to(V_1\oplus V_3,\omega_1-\omega_3)$ still carries an isotropic structure, thanks to the morphism $U\to L_{13}$. Since the isotropic structure on $L_{13}\to(V_1\oplus V_3,\omega_1-\omega_3)$ is non-degenerate,  
\[
U\oplus W[-1]\simeq L_{13}\simeq \hofib(V_1\oplus V_3\to L_{13}^*)\simeq \hofib(V_1\oplus V_3\to U^*)\oplus W^*\,,
\]
and thus $\hofib(V_1\oplus V_3\to U^*)\not\simeq U$ (because $W\neq0$). Therefore the isotropic structure on $U\to(V_1\oplus V_3,\omega_1-\omega_3)$ is not lagrangian. 

We therefore recover that the usual fiber product $U=\ker(L_{12}\oplus L_{23}\to V_2)$ is lagrangian if \textbf{and only if} the intersection is transverse. 
\hfill$\triangle$\end{Example}

\begin{Example}[odd shifted lagrangian intersection]\label{ex:1-shifted-lagrangian-intersection}
    Let $V_1=V_3=0$ again, and let $V_2$ be a $2$-term complex 
    \begin{equation*}
        E\overset{a}{\longrightarrow} F
    \end{equation*}
    with an $n$-shifted symplectic structure, for $n=-1-2d$ odd, as in \cref{subsec-2.2.2}. Here we assume that $\alpha=\omega_L^\flat: E\to F^*$ is an isomorphism.
    Let $L_{12}= \left(E \xrightarrow{b} B\right)$, a $2$-term lagrangian from \cref{subsec-2.3.4}, and $L_{23}= \left(0\to F\right)$. The latter admits a lagrangian structure on 
    \begin{equation*}
        \begin{pmatrix}
        0\\ id_F
        \end{pmatrix}
        :L_{23}\longrightarrow
        V_2,
    \end{equation*} 
    with respect to the obvious (meaning zero) isotropic structure, which is non-degenerate because 
    \begin{equation*}
    \left(
    \begin{tikzcd}
    0 \\
    F
    \arrow[from=1-1, to=2-1]
    \end{tikzcd}
    \right)
    \longrightarrow
    \hofib\left(
    \begin{tikzcd}[column sep=scriptsize]
	E & {F^*} \\
	F & 0
	\arrow["\alpha", from=1-1, to=1-2]
	\arrow["a"', from=1-1, to=2-1]
	\arrow[from=1-2, to=2-2]
	\arrow[from=2-1, to=2-2]
    \end{tikzcd}
    \right)
    =
    \begin{tikzcd}
	E \\
	{F\oplus F^*}
	\arrow["{(a,\alpha)}", from=1-1, to=2-1]
    \end{tikzcd}
    \end{equation*}
    is a quasi-isomorphism. Indeed, $\alpha$ is an isomorphism, so $\ker(a,\alpha) = 0$ and $\coker(a, \alpha)\cong F$. 
    Therefore, $L_{13}$ is $(n-1)$-shifted symplectic. 
    Recall 
    \begin{equation*}
    L_{13} = \hofib\left(
    \begin{tikzcd}
	E \\
	B
	\arrow["b"', from=1-1, to=2-1]
    \end{tikzcd}
    \oplus 
    \begin{tikzcd}
    0 \\
    F
    \arrow[from=1-1, to=2-1]
    \end{tikzcd}
    \longrightarrow
    \begin{tikzcd}
	E \\
	F
	\arrow["a", from=1-1, to=2-1]
    \end{tikzcd}
    \right),
    \end{equation*}
    which is
    \begin{equation}\label{eq-hofib}
    \begin{tikzcd}[column sep=scriptsize]
	{\text{degree}\quad  d} & & E & \\
	{\text{degree}\quad  d+1} & B & F & E \\
	{\text{degree}\quad  d+2} & & F &
	\arrow["b"', from=1-3, to=2-2]
	\arrow["{id_E}", from=1-3, to=2-4]
	\arrow["\oplus"{description}, draw=none, from=2-3, to=2-4]
	\arrow["f"', from=2-2, to=3-3]
	\arrow["\oplus"{description}, draw=none, from=2-2, to=2-3]
	\arrow["{id_F}", from=2-3, to=3-3]
	\arrow["a", from=2-4, to=3-3]
    \end{tikzcd}
    \end{equation}
    Here we see that the parts $E \xrightarrow{id_F} E$ and $F \xrightarrow{id_F} F$ are acyclic, so that \eqref{eq-hofib} is quasi-isomorphic to $B$ (sitting in degree $1+d$).
    Thus $B[-1-d]\simeq L_{13}$ is $(n-1)$-shifted symplectic. Moreover, keeping track of various identifications, one can prove that the $(n-1)$-shifted symplectic pairing is exactly the map $\eta:\wedge^2 B\to k$ defining the lagrangian structure on $L_{12}$. 
 \hfill$\triangle$\end{Example}


\section{Shifted symplectic (affine) derived schemes}\label{section3}

Recall the basic situation of a \emph{smooth affine algebraic variety} (of pure dimension) $X=\Spec(A)$: the algebra of functions $A:= \mathcal{O}(X)$ is 
finitely generated and the module of Kähler differentials $\Omega^1_A = \Omega^1(X)$ is a projective $A$-module of finite rank $\dim(X)$. 
This is equivalent to requiring that $A$ is finitely generated and the module of derivations
\begin{equation*}
    T_A:= \Der(A)=\VF(X)\,,
\end{equation*}
that is dual to $\Omega^1_A$, is a projective $A$-module of finite rank $\dim(X)$. 
We also say that $A$ is a \emph{smooth algebra}. 

The \emph{De Rham complex} of $X$ (or, of $A$) is 
\begin{equation*}
    \Sym_A(\Omega^1_A[-1]) = \bigoplus_{n\in\mathbb{N}}\Omega^n_A[-n],\quad\textrm{with}\quad \Omega^n_A=\wedge^n_A\Omega^1_A, 
\end{equation*}
with differential $d(a_0 da_1 \wedge \dots \wedge da_n) = da_0 \wedge da_1 \wedge \dots \wedge da_n$.

\begin{Definition}
    A \emph{symplectic structure} on a smooth affine algebraic variety $X=\Spec(A)$ is a $d$-closed $2$-form $\omega \in \Omega^2(X) = \Omega^2_A$, such that, viewed as a pairing $\wedge^2_AT_A\to A$, it makes $T_A$ a symplectic $A$-module.
\end{Definition}

Our main goal in this section is to extend this definition from smooth affine algebraic varieties to affine derived schemes, that is to the situation where $A$ is a \emph{connective} (meaning non-positively graded) \emph{commutative differential graded algebra}, or cdga for short. 

Before doing so, let us recall the following:
\begin{Idea}[Homotopical algebra]
Resolve problems \emph{before} they appear. Concretely this means that one must always take suitable resolutions before applying functors. 
\end{Idea}
The above idea lies at the origin of the yoga of derived functors, that is at the foundation of modern homological algebra \cite{CE}. It has been extended by Quillen \cite{QuillenHA} to the non additive setting. 
Earlier in these notes, we have already encountered an incarnation of this idea. 
\begin{Example}[Homotopy fiber product]\label{ex-hfp}
Let $V\xrightarrow{g} W$ be a linear map between vector spaces (or more generally a cochain map between complexes). The fiber product (or pullback) functor $-\times_WV$ is in general not well behaved. For instance, the expected dimension of the fiber product does not generally coincide with its actual dimension; but it is well-behaved whenever it is applied to a \emph{surjective} map $U\xrightarrow{f} W$. In Quillen's language, such a surjective map is a \emph{fibrant object} for a model structure on category of complexes over $W$. In case the map is not surjective, we consider the following \emph{fibrant resolution}  $f\to \tilde{f}$:  
\[\begin{tikzcd}
	\widetilde{U}:=U\oplus \hofib{(W\xrightarrow{id_W} W)} \\
	U & W
	\arrow[hook, "\simeq", from=2-1, to=1-1]
	\arrow["{(f,id_W,0)}=:\tilde{f}", from=1-1, to=2-2]
	\arrow["f", from=2-1, to=2-2]. 
\end{tikzcd}\]
Then, the \emph{homotopy} (also called \emph{right derived}) fiber product is 
\[
U\overset{h}{\times}_WV:=\widetilde{U}\times_WV\cong \hofib{(U\oplus V\overset{f-g}{\longrightarrow} W)}. 
\]
The general formalism ensures that different choices of resolutions give quasi-isomorphic results. 
\hfill$\triangle$\end{Example}

\subsection{A hint of derived geometry}

We denote the category of connective (i.e.~non positively graded) cdgas (commutative differential graded algebras) by $\cdga_k^{\le 0}$. 
\begin{Definition}
Let $f:C\to A$ be a morphism in $\cdga_k^{\le 0}$. A \emph{quasi-free resolution} of $f$ is a factorization $C\xrightarrow{\tilde{f}} \widetilde{A}\to A$ such that: 
\begin{enumerate}
    \item The morphism $\widetilde{A}\to A$ is a quasi-isomorphism. 
    \item Denoting by $(-)^{\not{\,\delta}}$ the ``underlying commutative graded algebra'' functor (that forgets the differential), we have that $\widetilde{A}^{\not{\,\delta}}\cong \Sym_{C^{\not{\,\delta}}}(V)$ as a commutative $C^{\not{\,\delta}}$-algebra, for some non-positively graded vector space $V$. 
\end{enumerate}
A quasi-free resolution of a cdga $A$ is a quasi-free resolution of the unit morphism $1:k\to A$. 
\end{Definition}
\begin{Remark}
    The categories $A\Mod$ and $\widetilde{A}\Mod$ are quasi-equivalent (quasi-equivalences of dg-categories are to equivalences of categories as quasi-isomorphisms are to isomorphisms) under the extension of scalars functor
    \begin{equation*}
    \begin{split}
        \widetilde{A}\Mod &\longrightarrow A\Mod\\
        M &\longmapsto A\otimes_{\widetilde{A}} M.
    \end{split}
    \end{equation*}
\end{Remark}
We are not going to detail the whole homotopy theory (i.e.~Quillen model structure) of $\cdga_k^{\le 0}$, as we only need the following: 
\begin{itemize}
    \item Weak equivalences of connective cdgas are quasi-isomorphisms of such. 
    \item Quasi-free resolutions always exist and are examples of cofibrant resolutions (i.e.~resolutions one uses to compute left derived functors in general, and homotopy pushouts more specifically). 
\end{itemize}
\begin{Remark}
Note that whenever the degree zero part $A^0$ of $A$ is finitely generated, one may just use a smooth resolution of $A$ instead of a quasi-free one: for a smooth resolution, we only require that $\widetilde{A}^0$ is a smooth algebra and $\widetilde{A}^{\not{\,\delta}}\cong\Sym_{\widetilde{A}^0}(P)$ with $P$ a negatively graded $\widetilde{A}^0$-module (which can be chosen to be projective). These specific models are exactly the (function dg-algebras of) the ``dg-manifolds'' from \cite[Definition 2.5.1]{Quot}. In differential geometry, this amounts to having $A^0=C^\infty(M)$ (with $M$ a smooth manifold) and $\widetilde{A}^{\not{\,\delta}}\cong\Gamma\big(M,\Sym(E)\big)$ with $E$ a negatively graded vector bundle (of finite rank in every degree) on $M$.   
\end{Remark}

\begin{Example}\label{ex-15}
    Consider the algebra $A=k[x]/x^2$, which is not smooth.  
    A quasi-free resolution of $A$ can be obtained by adding an extra generator $\xi$ in degree $-1$ and writing
    \begin{equation*}
        \widetilde{A} = k[x,\xi] = k[x]\xi \oplus k[x]= \left(k[x]\xi \overset{\delta}{\longrightarrow} k[x] \right)
    \end{equation*}
    with $\xi^2=0$ and differential $\delta$ given by $\delta(\xi) = x^2$. 
\hfill$\triangle$\end{Example}

\subsubsection{Tangent and cotangent complexes}

\begin{Definition}
    Consider $A\in \cdga_k^{\le 0}$, and a quasi-free resolution $\widetilde{A}$. We define the \emph{tangent} and \emph{cotangent complexes}, respectively, as 
    \begin{equation*}
        \Tang_A := T_{\widetilde{A}}, 
        \quad
        \Cotang_A := \Omega^1_{\widetilde{A}} 
    \end{equation*}
in $\widetilde{A}\Mod\simeq A\Mod$. 
\end{Definition}

\begin{Example}\label{ex-fat-point}
    Using the quasi-free resolution from \cref{ex-15}, the cotangent complex of $A=k[x]/x^2$ is 
    \begin{equation*}
        \Cotang_A = k[x,\xi]d\xi \oplus k[x,\xi]dx,
    \end{equation*}
    with $\deg(d\xi)=-1$, and differential given by 
    \begin{equation*}
        \delta(\xi) = x^2\qquad \textrm{and} \qquad \delta(d\xi) = d(\delta\xi) = 2xdx.
    \end{equation*}
    
    The tangent complex of $A$ is 
    \begin{equation*}
        \Tang_A = k[x,\xi]\frac{\partial}{\partial x} \oplus k[x,\xi]\frac{\partial}{\partial \xi}, 
    \end{equation*}
    with $\deg\left(\frac{\partial}{\partial \xi}\right)= 1$, and 
    \begin{equation*}
        \langle \delta \left(\frac{\partial}{\partial x}\right), d\xi \rangle = \langle \frac{\partial}{\partial x}, \delta d \xi \rangle = \langle \frac{\partial}{\partial x}, 2xdx \rangle= 2x, 
    \end{equation*}
    i.e. $\delta \left(\frac{\partial}{\partial x}\right) = 2x \frac{\partial}{\partial \xi}$.
    We observe that $\Tang_A$ is a $(-1)$-shifted symplectic $\widetilde{A}$-module, with respect to the unique skew-symmetric pairing such that 
    \begin{equation*}
        \omega\left(\frac{\partial}{\partial x} \wedge \frac{\partial}{\partial \xi}\right) = 1.
    \end{equation*}
\hfill$\triangle$\end{Example}

\subsubsection{Derived fiber products}

Following \cref{ex-hfp}, we are going to define the relative (left) derived tensor product (or homotopy pushout) in $\cdga_k^{\le 0}$; in other words, for $C\xrightarrow{g} B$ a morphism in $\cdga_k^{\le 0}$, we are going to derive the functor $-\otimes_C B$. 
\begin{Definition}
    Let $C\xrightarrow{f} A$ be another morphism in $\cdga_k^{\le 0}$. The relative (left\footnote{We will not explain in these notes the difference between left and right derived functors. }) \emph{derived tensor product} of $A$ and $B$ over $C$ is 
    \begin{equation*}
        A\overset{\mathbf{L}}{\otimes_C} B:= \widetilde{A} \otimes_C B,
    \end{equation*}
    where $C\xrightarrow{\tilde{f}} \widetilde{A}$ is a quasi-free resolution of $f$. 
\end{Definition}
The general formalism of model categories ensures that different choices of resolutions give quasi-isomorphic results. 
\begin{Remark}
    Note that the relative tensor product (i.e.~pushout) of commutative algebras is the algebraic incarnation, at the level of functions, of the fiber product (i.e.~pullback) of affine schemes. 
    Hence the relative derived tensor product defines a \emph{derived/homotopy fiber product} (often abusively called \emph{derived intersection}) of affine derived schemes. 
    Geometrically, if $A=\mathcal O(X)$ and $C=\mathcal O(Z)$, then the quasi-free replacement $C\overset{\tilde{f}}{\longrightarrow} \widetilde{A}\overset{\simeq}{\longrightarrow}A$ gives a factorization $X\overset{\simeq}{\longrightarrow} \widetilde{X}\longrightarrow Z$, where $\tilde{X}\to Z$ encodes a submersion over $Z$ (algebraic geometers would say a \emph{smooth morphism} to $Z$) together with equations (encoded within the differential) that cut out $X$. This very idea dates back to Ciocan-Fontanine--Kapranov \cite[Theorem 2.7.6 \& Subsection 2.8]{Quot}. 
\end{Remark}

A nice feature of derived tensor products is that the tangent complex sends derived tensor products of finite type connective cdgas\footnote{The more fundamental result is that the cotangent complex sends homotopy pushouts (i.e.~derived tensor products) to homotopy pushouts. One then needs to dualize to get the analogous result for tangent complexes, and duality sends pushouts to pullbacks under some finiteness assumption. } to derived fiber products. 
This is built in the construction: the ordinary tangent indeed sends pullbacks along smooth (i.e.~submersive) morphisms
to pullbacks. Being more precise, if $D:= A \overset{\mathbf{L}}{\otimes_C} B$ then 
\begin{equation*}
    \Tang_D \simeq D\otimes_A \Tang_A\overset{h}{\underset{D \otimes_C \Tang_C}{\times}}D\otimes_B \Tang_B.
\end{equation*}

\begin{Example}\label{ex-long}
    Consider $X=\mathbb{A}^1 \hookrightarrow \mathbb{A}^2= Z$ the affine line embedded into the affine plane as $\{y=0\}$. We would like to compute the derived self-intersection of $X$ into $Z$. Algebraically, on functions, we have 
    \begin{equation*}
    \begin{split}
        A:=\mathcal O(X)=k[x] &\longleftarrow k[x,y]=\mathcal O(Z)=:C \\
        0 &\longleftarrow\!\shortmid y.
    \end{split}
    \end{equation*}
    The cdga of functions on the derived self-intersection of $\mathbb{A}^1$ in $\mathbb{A}^2$ is computed by the derived tensor product $A \overset{\mathbf{L}}{\otimes_C} A$:
    \begin{enumerate}
        \item One first resolves $C \to A$ by considering 
        \begin{equation*}
            C\longrightarrow \widetilde{A}=k[x,y,\xi] \longrightarrow A,
        \end{equation*}
        with $\deg(y)=0$, $\deg(\xi)=-1$, and $\delta \xi = y$. 
        \item Then the derived tensor product is 
        \[
        A \overset{\mathbf{L}}{\otimes_C} A=\widetilde{A}\otimes_C A \cong k[x,\xi]=k[x]\otimes k[\xi]
        \]
        (with $\delta(\xi)=0$ now), which corresponds to the space $\mathbb{A}^1 \times \mathbb{A}^1[-1]$, where $\mathbb{A}^1[-1]$ is the odd affine line. 
    \end{enumerate}
The tangent complex of the derived self-intersection is
\begin{equation*}
            \Tang_{k[x,\xi]} = k[x,\xi] \frac{\partial}{\partial x} \oplus k[x,\xi] \frac{\partial}{\partial \xi}.
\end{equation*}
with $\deg(\frac{\partial}{\partial \xi})=1$ and zero differential. We now compare this with 
\begin{align*}
& k[x,\xi]\otimes_A \Tang_A\overset{h}{\underset{k[x,\xi] \otimes_C \Tang_C}{\times}}k[x,\xi]\otimes_A \Tang_A \\
\cong & \hofib\left((k[x,\xi]\otimes_{\widetilde{A}} T_{\widetilde{A}})^{\oplus2}  \overset{(t,-t)}\longrightarrow k[x,\xi]\otimes_{C} T_{C} \right).
\end{align*}
Here we have that: 
\begin{enumerate}
\item The differential graded $k[x,\xi]$-module $k[x,\xi]\otimes_{\widetilde{A}} T_{\widetilde{A}}$ is freely generated by $\frac{\partial}{\partial x}$ (in degree $0$), $\frac{\partial}{\partial y}$ (in degree $0$) and $\frac{\partial}{\partial \xi}$ (in degree $1$), with $\delta(\frac{\partial}{\partial y})=\frac{\partial}{\partial \xi}$. 
\item The differential graded $k[x,\xi]$-module $k[x,\xi]\otimes_{C} T_{C}$ is freely generated in degree $0$ by $\frac{\partial}{\partial x}$ and $\frac{\partial}{\partial y}$. 
\item The morphism (of differential graded $k[x,\xi]$-modules) $t$ sends $\frac{\partial}{\partial \xi}$ to $0$ and is the identity on the other generators. 
\end{enumerate}
Hence $\hofib\left((k[x,\xi]\otimes_{\widetilde{A}} T_{\widetilde{A}})^{\oplus2}  \overset{(t,-t)}\longrightarrow k[x,\xi]\otimes_{C} T_{C} \right)$ is the free differential graded $k[x,\xi]$-module generated by the following $2$-term complex: 
\[\begin{tikzcd}
	{\text{degree 0}} & \textcolor{rgb,255:red,163;green,41;blue,41}{{k\frac{\partial}{\partial x}}} & \textcolor{rgb,255:red,163;green,41;blue,41}{{k\frac{\partial}{\partial y}}} & \textcolor{rgb,255:red,24;green,98;blue,24}{{k\frac{\partial}{\partial x}}} & \textcolor{rgb,255:red,24;green,98;blue,24}{{k\frac{\partial}{\partial y}}} \\
	{\text{degree 1}} & \textcolor{rgb,255:red,163;green,41;blue,41}{{k\frac{\partial}{\partial\xi}}} & \textcolor{rgb,255:red,33;green,84;blue,171}{{k\frac{\partial}{\partial x}}} & \textcolor{rgb,255:red,33;green,84;blue,171}{{k\frac{\partial}{\partial y}}} & \textcolor{rgb,255:red,24;green,98;blue,24}{{k\frac{\partial}{\partial\xi}}}
	\arrow["\oplus"{description}, draw=none, from=1-2, to=1-3]
	\arrow[from=1-2, to=2-3]
	\arrow["\oplus"{description}, draw=none, from=1-3, to=1-4]
	\arrow[color={rgb,255:red,163;green,41;blue,41}, from=1-3, to=2-2]
	\arrow[from=1-3, to=2-4]
	\arrow[from=1-4, to=2-3]
	\arrow[from=1-5, to=2-4]
	\arrow[color={rgb,255:red,24;green,98;blue,24}, from=1-5, to=2-5]
	\arrow["\oplus"{description}, draw=none, from=2-2, to=2-3]
	\arrow["\oplus"{description}, draw=none, from=2-3, to=2-4]
\end{tikzcd}\]
where the red (resp.~green) part corresponds to the first (resp.~second) copy of $k[x,\xi]\otimes_{\widetilde{A}} T_{\widetilde{A}}$, the blue part corresponds to $k[x,\xi]\otimes_{C} T_{C}$, and the remaining arrows describe the map $(t,-t)$. 
Observe that the above $2$-term complex projects to $k\frac{\partial}{\partial x}\oplus k\frac{\partial}{\partial \xi}$, with kernel 
\[\begin{tikzcd}
	{\text{degree 0}} && \textcolor{rgb,255:red,163;green,41;blue,41}{{k\frac{\partial}{\partial y}}} & \textcolor{rgb,255:red,24;green,98;blue,24}{{k\frac{\partial}{\partial x}}} & \textcolor{rgb,255:red,24;green,98;blue,24}{{k\frac{\partial}{\partial y}}} \\
	{\text{degree 1}} && \textcolor{rgb,255:red,33;green,84;blue,171}{{k\frac{\partial}{\partial x}}} & \textcolor{rgb,255:red,33;green,84;blue,171}{{k\frac{\partial}{\partial y}}} & \textcolor{rgb,255:red,24;green,98;blue,24}{{k\frac{\partial}{\partial\xi}}}
	\arrow["\oplus"{description}, draw=none, from=1-3, to=1-4]
	\arrow[from=1-3, to=2-4]
	\arrow[from=1-4, to=2-3]
	\arrow[from=1-5, to=2-4]
	\arrow[draw={rgb,255:red,24;green,98;blue,24}, from=1-5, to=2-5]
	\arrow["\oplus"{description}, draw=none, from=2-3, to=2-4], 
\end{tikzcd}\]
which is acyclic, ensuring that the projection mentioned above is a quasi-isomorphism. We therefore get that 
\[
k[x,\xi]\otimes_A \Tang_A\overset{h}{\underset{k[x,\xi] \otimes_C \Tang_C}{\times}}k[x,\xi]\otimes_A \Tang_A\simeq \Tang_{k[x,\xi]}
\]
as expected. 
\hfill$\triangle$\end{Example}

\subsection{De Rham complex and shifted symplectic structures}

\begin{Definition}[Pantev--Toën--Vaquié--Vezzosi  \cite{PTVV}]
    Let $(A,\delta) \in \cdga_k^{\le 0}$. The \emph{de Rham complex} of $A$ is 
    \begin{equation*}
        DR^\bullet (A) := \Sym_{\widetilde{A}}^\bullet\left(\Omega_{\widetilde{A}}^1[-1]\right). 
    \end{equation*}
    This is a graded mixed complex\footnote{Up to a shift in the grading conventions, this is the same as a bicomplex. }, as it has two commuting differentials:
    \begin{enumerate}
        \item The internal differential $\delta$, which preserves the symmetric weight. 
        \item The de Rham differential $d$, which increases the symmetric weight by 1. 
    \end{enumerate}
\end{Definition}

\begin{Definition}[Pantev--Toën--Vaquié--Vezzosi \cite{PTVV}]\label{def:n-shifted-forms}
\hfill
\begin{enumerate}
    \item A \emph{$2$-form of degree $n$} on $A$ is an $(n+2)$-cocycle in 
    \begin{equation*}
        \Sym_{\widetilde{A}}^2\left(\Omega_{\widetilde{A}}^1[-1]\right) \cong \wedge^2_{\widetilde{A}}\left(\Omega_{\widetilde{A}}^1 \right)[-2],
    \end{equation*}
    with respect to $\delta$.
     Such a $2$-form induces a cochain map (i.e. a pairing)
    \begin{equation*}
        \wedge^2_{\widetilde{A}} T_{\widetilde{A}} \longrightarrow \widetilde{A}[n].
    \end{equation*}
    The $2$-form is said to be \emph{non-degenerate} if this pairing is. 

    \item A \emph{closed $2$-form of degree $n$}  is an $(n+2)$-cocycle in the total complex
    \begin{equation*}
        \left(\prod_{i \ge 0} \Sym_{\widetilde{A}}^{2+i}\left(\Omega_{\widetilde{A}}^1 [-1]\right),\delta+d\right).
    \end{equation*}
    Explicitly, this is a series $\omega_0 + \omega_1 + \dots$ with $\omega_i$ of weight $2+i$, such that 
    \begin{equation*}
        \delta\omega_0=0, \text{ and } \forall i\ge 0, \; d\omega_i + \delta \omega_{i+1} = 0.
    \end{equation*}
    In particular, $\omega_0$ is a $2$-form of degree $n$ and $\omega_{1}$ is a homotopy between $0$ and $d\omega_1$. The other $\omega_i$'s are higher coherent homotopies. 

    \item An \emph{$n$-shifted symplectic structure} on $A$ (or $X= \Spec(A)$) is a closed $2$-form of degree $n$ such that $\omega_0$ is non-degenerate. 
\end{enumerate}
\end{Definition}

\begin{Example}
    Assume that $A$ carries a $0$-shifted symplectic structure. 
    Then $\omega_0^\flat:\Tang_A \to \Cotang_A$ is a quasi-isomorphism. 
    Because $\Tang_A$ is concentrated in non-negative degrees, and $\Cotang_A$ is concentrated in non-positive degrees, we have a sequence of quasi-isomorphisms 
    \begin{equation*}
        H^0(\Tang_A) 
        \xrightarrow{q. iso.}
        \Tang_A 
        \xrightarrow{q. iso.}
        \Cotang_A
        \xrightarrow{q. iso.}
        H^0(\Cotang_A). 
    \end{equation*}
    Therefore, up to quasi-isomorphism, $A$ is concentrated in degree $0$ and it is smooth. Thus it is (equivalent to) an honest symplectic scheme. 
\hfill$\triangle$\end{Example}

\begin{Example}
    Assume $A$ carries a $(-1)$-shifted symplectic structure. 
    Then 
    \begin{equation*}
        \omega_0^\flat: \;
        \Tang_A \longrightarrow
        \Cotang_A[-1]
    \end{equation*}
    is a quasi-isomorphism. 
    Thus $\Tang_A$ has cohomology in degrees 0 and 1, of finite rank, by the same argument as before. 
    In this case $A$ is called \emph{quasi-smooth}, as the defect of smoothness is controlled by a single vector bundle: the obstruction bundle. 
\hfill$\triangle$
\end{Example}

\begin{Remark}
    For obvious degree reasons, there are no $n$-shifted symplectic structures on $A$ for $n>0$ (apart from the zero $n$-shifted symplectic structure on $k$). 
\end{Remark}

\begin{Remark}
    Locally, $n$-shifted symplectic structures on $A$ have strict normal forms. That is, for every $k$-point there exists a Zariksi open neighborhood $U$ and a quasi-free model $\widetilde{A}\simeq A_U$ for which $d\omega_0=0$, $\omega_i=0$ for all $i\ge 1$, and $\omega_0^\flat$ is an isomorphism. 
    This was shown by Brav--Bussi--Joyce \cite{BBJ}.
\end{Remark}

\subsection{Lagrangian morphisms}

From now on, we will think of connective cdgas in geometric terms. Formally, this means that we consider the category of \emph{affine derived schemes}, which is the opposite of the homotopy category\footnote{The homotopy category is obtained by formally inverting quasi-isomorphisms. Since $\cdga_k^{\le 0}$ carries a Quillen model structure, it can be described in more concrete terms: object are quasi-free cdgas, and hom sets are quotients of hom sets in $\cdga_k^{\le 0}$ by a homotopy equivalence relation. } of $\cdga_k^{\le 0}$. 
\begin{Definition}\label{definition:lagrangian-affine}
    Let $f: Y \to X$ be a morphism of affine derived schemes represented by $f^*: A \to B$ in $\cdga_k^{\le 0}$. 
    Assume that $X$ (i.e. $A$) is equipped with an $n$-shifted symplectic structure $\omega$. 
    \begin{enumerate}
    \item An \emph{isotropic structure} on $f$ with respect to $\omega$ is a homotopy $\eta$ between $f^*\omega$ and $0$. 
    Concretely, 
    \begin{equation*}
        \eta = \eta_0 + \eta_1 + \eta_2 + \dots,
    \end{equation*}
    with $\eta_i$ of weight $2+i$ and $f^*\omega = (d+\delta)(\eta)$. 
    In particular $f^*\omega_0=\delta\eta_0$, meaning that the $n$-shifted pairing \[
    f^*\omega_0:\wedge_B\mathbb{T}_B\to B\otimes_A(\wedge^2_A\Tang_A)\xrightarrow{id_B\otimes\omega_0}B\otimes_AA[n]=B[n]
    \]
    on $\Tang_B$ is homotopic to zero, \emph{via} $\eta_0$. 
    \item A \emph{lagrangian structure} on $f$ with respect to $\omega$ is an isotropic structure $\eta$ such that $\eta_0$ is non-degenerate in the sense of  \cref{def-nd-homotopy}. 
    \end{enumerate}
\end{Definition}

This definition of lagrangian morphisms for affine derived schemes presents many of the same features as in the linear setting, such as:
\begin{enumerate}
    \item Lagrangian structures on $X \to pt:=\Spec(k)$ with respect to the null $n$-shifted symplectic structure on $pt$ are $(n-1)$-shifted symplectic structures on $X$ \cite[Example 2.3]{Cal}. 
    \item Lagrangian correspondences compose well, provided one uses \emph{derived} fiber products \cite[Theorem 4.4]{Cal}.
    \item Genuine smooth lagrangian subschemes are examples of lagrangian morphisms \cite[proof of Corollary 2.10]{PTVV}. 
\end{enumerate}
As a consequence, the derived fiber product $X\overset{h}{\times}_Z Y$ of two lagrangian morphisms $X\rightarrow Z\leftarrow Y$ of affine derived schemes towards the same $n$-shifted symplectic affine derived scheme $Z$ is $(n-1)$-shifted symplectic. In particular, derived lagrangian intersections are $(-1)$-shifted symplectic, a fact that we are going to illustrate in in the remainder of this section.  

\subsubsection{Derived zero locus of a closed $1$-form}

Let $X= \Spec(A)$ be a smooth affine algebraic variety. Its cotangent bundle $T^*X$ is symplectic with respect to the canonical form. 
The graph of any closed $1$-form $\lambda: X \to T^*X$ (such as for instance the zero section $0: X\to T^*X$) is lagrangian. One can then consider the \emph{(right) derived zero locus} of a closed $1$-form $\lambda$, 
\begin{equation*}
    \mathbf{R} Z(\lambda) := X \overset{h}{\underset{0, T^*X, \lambda}{\times}}X.
\end{equation*}
This is $(-1)$-shifted symplectic, as it is a derived intersection of lagrangian subvarieties. 
Whenever $\lambda = df$, $\mathbf{R} Z(\lambda) =: \mathbf{R} \mathrm{Crit}(f)$ is the \emph{derived critical locus of $f$}. 

\begin{Example}
Let $X=\Spec(k[x])$ and $\lambda = \alpha(x)dx$. In this case, 
\begin{equation*}
    T^*X \cong\mathbb{A}^2=\Spec(k[x,y]), 
    \text{ with }
    \omega= dx \wedge dy = \omega_0.
\end{equation*}
As in \cref{ex-long} we use the quasi-free resolution 
\[
k[x,y]\longrightarrow k[x,y,\xi]\overset{\simeq}{\longrightarrow}k[x], 
\]
with $\delta(\xi)=y$, and get 
\[
\mathcal O\big(\mathbf{R} Z(\lambda)\big)=k[x,y,\xi]\otimes_{k[x,y]}k[x]\cong \big(k[x,\xi],\delta:\xi\mapsto \alpha(x)\big),
\]
because the morphism $k[x,y]\to k[x]$ representing the section $\lambda$ sends $y$ to $\alpha(x)$. 
The differential can thus be written as $\delta=\alpha\frac{\partial}{\partial\xi}$. 

In order to compute the $(-1)$-shifted symplectic structure on the derived intersection, we first have to understand what is the homotopy between $\omega$ and $0$ in the quasi-free resolution $k[x,y,\xi]$. One easily sees that it is $\eta_0=\eta=dx\wedge d\xi$:  
\begin{equation*}
    d\eta_0 = 0, \text{ and }
    \delta \eta_0 = 
    \delta d x \wedge d\xi - dx \wedge \delta d\xi = dx \wedge dy = \omega.
\end{equation*}
Therefore, the $(-1)$-shifted symplectic structure on $(k[x,\xi], \alpha\frac{\partial}{\partial\xi})$ is $dx\wedge d\xi$. 
\hfill$\triangle$\end{Example}
The computation from the above example can be adapted to the more general situation and one can prove that the cdga of functions on $\mathbf{R} Z(\lambda)$ is given by 
\begin{equation*}
    \mathcal O\big(\mathbf{R} Z(\lambda)\big) \simeq \big(\Sym_A (T_A[1]), \iota_\lambda\big),
\end{equation*}
where $\iota_\lambda$ is the contraction with respect to $\lambda$. The $(-1)$-symplectic structure is the canonical odd symplectic structure on $\Sym_A (T_A[1])=\mathcal O(T^*[-1]X)$ known to Poisson geometers (one can check that it is closed independantly for the internal differential and for the de Rham differential). 
\begin{Example}
    If $\lambda = df = x^2dx$, that is, $f= \frac{1}{3} x^3$, then 
    \begin{equation*}
        \mathcal O\big(\mathbf{R} Z(\lambda)\big) = \mathcal O\big(\mathbf{R} \mathrm{Crit}(\frac{1}{3} x^3)\big)\simeq (k[x,\xi],x^2\frac{\partial}{\partial\xi})\simeq k[x]/x^2 
    \end{equation*}
carries a $(-1)$-shifted symplectic structure. As a consequence, we recover a fact already noticed in \cref{ex-fat-point}: the tangent complex of $k[x]/x^2$ carries a (linear) $(-1)$-shifted symplectic structure. 
\hfill$\triangle$\end{Example}


\section{Shifted symplectic 1-stacks}\label{section4}

\subsection{Higher and derived stacks}

Roughly speaking, higher stacks are sheaves with values in $\infty$-groupoids, in which case, the gluing axiom is only intended to hold up to homotopy. Depending on the setting, these are sheaves on different categories. 
\begin{itemize}
    \item In differential geometry, one considers sheaves over the Euclidean site, consisting of smooth manifolds with the topology of local diffeomorphisms. 
    \item In algebraic geometry, one often considers sheaves over the étale site, i.e. affine schemes with the topology of étale morphisms. 
    \item In derived algebraic geometry, one often considers the derived étale site, which consists of affine derived schemes with étale morphisms. 
\end{itemize}

These are often too general to be geometrically tractable. Therefore one usually considers \emph{geometric stacks}. These are the ones that can be obtained from representable objects (building blocks) by iterated smooth groupoid quotients.

\begin{Definition}
A \emph{smooth groupoid} is a simplicial object $X_\bullet$ such that 
\begin{enumerate}
    \item The face maps (also called \emph{source} and \emph{target}) $X_1 \rightrightarrows X_0$ are smooth morphisms;  
    \item The canonical projection to the space of $n$-tuples of composable 1-simplices
\begin{equation*}
    X_n \longrightarrow X_1 \underset{X_0}{\times} \dots  \underset{X_0}{\times} X_1 \text{ ($n$ times)}
\end{equation*}
is an equivalence for any $n$; 
    \item The canonical projection $X_2\to X_1\underset{X_0}{\times}X_1$ to the space of pairs of $1$-simplices with same source (resp.~target) is an equivalence. 
\end{enumerate}
\end{Definition}

If $X_n$ is a manifold (or a smooth affine algebraic variety) for every $n$ then this defines a \emph{Lie groupoid}. Indeed, conditions 2 and 3 are equivalent to requiring that $X_\bullet$ is the nerve of a groupoid. 
In the following discussion, we will only deal with those, i.e. with \emph{underived smooth 1-stacks}.

\begin{Remark}
In differential geometric terms, the first condition says that the source and target maps are submersive with finite-dimensional fibers. Smoothness implies in particular that that 
\begin{equation*}
    \hofib(\Tang_{X_1} \to \Tang_{X_0})
\end{equation*}
is perfect and concentrated in degree 0, and that the fiber product $X_1 \underset{X_0}{\times} \dots  \underset{X_0}{\times} X_1$ appearing in the second condition coincides with the derived fiber product. 
\end{Remark}

\medskip

The quotient of a groupoid $X_\bullet$ is defined as
\begin{equation*}
    |X_\bullet|=\underset{[n]\in\Delta^{op}}{\colim}(X_n).
\end{equation*}
\begin{Notation}
Whenever $X_\bullet$ is an underived smooth groupoid, one often writes $[X_0/X_1]$ for $|X_\bullet|$. Furthermore, if it is the action groupoid of a group $G$ acting on $X_0$ one would rather write $[X_0/G]$. Finally, whenever $X_0=*$, and therefore $X_1 = G$ is a group, we write $BG := [*/G]$. 
\end{Notation}

\begin{Remark}
    For honest Lie groupoids $X_\bullet$ and $Y_\bullet$, we have that the homotopy colimits $|X_\bullet|$ and $|Y_\bullet|$ are equivalent as stacks if and only if the groupoids are Morita equivalent (see e.g.~\cite{Metzler}). 
\end{Remark}

\subsection{De Rham complex and shifted symplectic structures}

The Yoneda lemma tells us that any (pre-)sheaf is the colimit of representables mapping to it (beware that for presheaves of $\infty$-groupoids, one should consider the \textit{homotopy} colimit). 
If $\mathcal{X}$ is a stack, then 
\begin{equation*}
    \mathcal{X} = \underset{\Spec A \to \mathcal{X}}{\colim}( \Spec A ).
\end{equation*}

\begin{Definition}
    The \emph{de Rham complex} of $\mathcal{X}$ is the limit 
    \begin{equation*}
        DR^\bullet(\mathcal{X}):= \underset{\Spec A \to \mathcal{X}}{\holim}(DR^\bullet(A))
    \end{equation*}
    in the category of graded mixed complexes.
\end{Definition}

Unfortunately, this definition is very convenient for theoretical purposes but not very practical for computations. Nevertheless, when $\mathcal{X}$ is nice enough (i.e. it admits a tangent and cotangent complex), then we have a quasi-isomorphism of graded complexes (thanks to \cite{Calaque-Safronov})
\begin{equation*}
    DR^\bullet(\mathcal{X})^{\not{d}} \simeq \Gamma(X, \Sym^\bullet_{\mathcal{O}_X}(\Cotang_X[-1])).
\end{equation*}
Therefore, the definitions of (closed) $2$-forms and symplectic structures from \cref{def:n-shifted-forms} apply verbatim in this setting, provided one has an explicit description of the de Rham differential on the graded complex of forms. We now describe them concretely for the case of (underived) geometric stacks.

\subsubsection{Concrete description for (underived) geometric stacks}

Assume that $\mathcal{X}=|X_\bullet|$, for $X_\bullet$ a smooth groupoid.\footnote{The same works for $X_\bullet$ a smooth $n$-groupoid.} 
Then 
\begin{equation*}
    DR(\mathcal{X})\simeq\underset{[n]\in\Delta}{\holim}DR(X_n) \simeq DR(X_0) 
    \xrightarrow{\widecheck{\partial}}
    DR(X_1)[-1] \xrightarrow{\widecheck{\partial}} 
    DR(X_2)[-2] 
    \rightarrow \dots
\end{equation*}
The internal differential is now $\delta + \widecheck{\partial}$, and the de Rham differential is still $d=d_{dR}$. 

\begin{Remark}
The first equivalence in the above equation comes from a smooth descent result (see \cite{PTVV}), while the second one is a standard computation of homotopy colimits of cosimplicial diagrams.  
\end{Remark}

For underived smooth stacks we have the following assumptions: $X_n$ is underived, not stacky and smooth. Therefore the complexes $\Tang_{X_n} \simeq T_{X_n}$ and $\Cotang_{X_n} \simeq \Omega^1_{X_n}$ are concentrated in degree 0. In this case we also have that $\delta = 0$, so we only have $\widecheck{\partial}$ and $d$. Thus we can draw the graded mixed de Rham complex as
\[\begin{tikzcd}
	& {\text{weight 0}} & {\text{weight 1}} & {\text{weight 2}} \\
	{\text{total degree 0}} & {\mathcal{O}(X_0)} \\
	{\text{total degree 1}} & {\mathcal{O}(X_1)} & {\Omega^1(X_0)} \\
	{\text{total degree 2}} & {\mathcal{O}(X_2)} & {\Omega^1(X_1)} & {\Omega^2(X_0)} \\
	{\text{total degree 3}} & {\mathcal{O}(X_3)} & {\Omega^1(X_2)} & {\Omega^2(X_1)} \\
	& \vdots & {\Omega^1(X_3)} & {\Omega^2(X_2)} \\
	&& \vdots & {\Omega^2(X_3)} \\
	&&& \vdots
	\arrow["{\widecheck{\partial}}"', from=2-2, to=3-2]
	\arrow["d", from=2-2, to=3-3]
	\arrow["{\widecheck{\partial}}"', from=3-2, to=4-2]
	\arrow["d", from=3-2, to=4-3]
	\arrow["{\widecheck{\partial}}"', from=3-3, to=4-3]
	\arrow["d", from=3-3, to=4-4]
	\arrow["{\widecheck{\partial}}"', from=4-2, to=5-2]
	\arrow["d", from=4-2, to=5-3]
	\arrow["{\widecheck{\partial}}"', from=4-3, to=5-3]
	\arrow["d", from=4-3, to=5-4]
	\arrow["{\widecheck{\partial}}"', from=4-4, to=5-4]
	\arrow[from=5-2, to=6-2]
	\arrow["d", from=5-2, to=6-3]
	\arrow["{\widecheck{\partial}}"', from=5-3, to=6-3]
	\arrow["d", from=5-3, to=6-4]
	\arrow["{\widecheck{\partial}}"', from=5-4, to=6-4]
	\arrow[from=6-3, to=7-3]
	\arrow["d", from=6-3, to=7-4]
	\arrow["{\widecheck{\partial}}"', from=6-4, to=7-4]
	\arrow[from=7-4, to=8-4]
\end{tikzcd}\]
Then, a $2$-form of degree $n$ on $\mathcal{X}=|X_\bullet|$ is a $2$-form $\omega \in \Omega^2(X_n)$ such that $\widecheck{\partial}\omega = 0$.
Meanwhile, a \emph{closed} $2$-form of degree $n$ on $\mathcal{X}=|X_\bullet|$ is a sum 
\begin{equation*}
    \omega= \omega_0 + \omega_1 + \omega_2 + \cdots +\omega_n
\end{equation*}
with $\omega_i \in \Omega^{2+i}(X_{n-i})$ and $(\widecheck{\partial} + d) (\omega) = 0$.

\subsubsection{The non-degeneracy condition in concrete terms}

We borrow the same assumptions as above: $X_n$ is underived, not stacky and smooth. Observe that we have a canonical quotient morphism 
\begin{equation*}
    p: X_0 \longrightarrow |X_\bullet|=\mathcal{X}.
\end{equation*}
One can check that the pullback $p^*$ is conservative, that is, 
$E \to F$ is a quasi-isomorphism of (quasi-coherent) sheaves on the quotient stack $\mathcal{X}$ if and only if $p^*E \to p^*F$ is a quasi-isomorphism\footnote{This is just saying that an equivariant morphism between $G$-sheaves on $X_0$ is a quasi-isomorphism if and only if it is so as a plain morphism of sheaves.}. 

\begin{Example}
    Consider $\mathcal{X} = BG = [*/G]$, for a Lie (or affine algebraic) group $G$. Then we have $p: * \to [*/G]$. Sheaves on $BG$ are cochain complexes of $G$-representations, and $p^*$ is the forgetful map, i.e. the map that forgets the $G$-action. A $G$-equivariant map between $G$-representations $E\to F$ is a quasi-isomorphism if and only it is a quasi-isomorphism of complexes. \hfill$\triangle$
\end{Example}

With this, the non-degeneracy condition can be checked after applying $p^*$. Several observations can be made at this point:
\begin{enumerate}
    \item By the above observation (that $p^*$ is conservative), it is sufficient to pull everything back to $X_0$ to check non-degeneracy. 
    \item The object $\Cotang_{X_\bullet|_{X_0}}\simeq \Omega^1_{X_\bullet|_{X_0}}$ is a cosimplicial sheaf on $X_0$, and 
    \begin{equation*}
        p^*\Cotang_{\mathcal{X}} \simeq
        \underset{[n]\in\Delta}{\holim}(\Omega^1_{X_n|_{X_0}})
        \simeq \Omega^1_{X_0}
        \xrightarrow{\widecheck{\partial}} 
        \Omega^1_{X_1|_{X_0}}[-1]
        \xrightarrow{\widecheck{\partial}} 
        \Omega^1_{X_2|_{X_0}}[-2]
        \to \dots,
    \end{equation*}
    starting with $\Omega^1(X_0)$ in degree 0 and increasing in degree from left to right. 
    Dually, 
    \begin{equation*}
        p^*\Tang_{\mathcal{X}} \simeq \underset{[n]\in\Delta^{op}}{\colim} (T_{X_\bullet|_{X_0}}
        \simeq \left(\dots \to 
        T_{X_2|_{X_0}}[2]\to T_{X_1|_{X_0}}[1] \to 
        T_{X_0}\right), 
    \end{equation*}
    where $T_{X_1|_{X_0}}$ sits in degree $-1$ and $T_{X_0}$ in degree 0.
    \item Since $X_\bullet$ is a groupoid, $X_n \simeq X_1 \times_{X_0} \dots \times_{X_0} X_1$ and thus 
    \begin{equation*}
        T_{X_n|_{X_0}} \cong T_{X_1|_{X_0}} \underset{T_{X_0}}{\times} \dots \underset{T_{X_0}}{\times} T_{X_1|_{X_0}}, 
        \text{ and }
        \Omega^1_{X_n|_{X_0}} \cong \Omega^1_{X_1|_{X_0}} \underset{X_0}{\oplus} \dots \underset{X_0}{\oplus} \Omega^1_{X_1|_{X_0}}.
    \end{equation*}
    More importantly, $T_{X_\bullet|_{X_0}}$ is the nerve of a groupoid object in the category of sheaves on $X_0$: $T_{X_1|_{X_0}} \overset{s_*}{\underset{t_*}{\rightrightarrows}} T_{X_0}$. It is a general result that for a groupoid object $B \overset{f}{\underset{g}{\rightrightarrows}} A$ in a stable $\infty$-category, the homotopy colimit of its nerve is $\hofib(B\overset{g-f}{\longrightarrow}A)$; therefore $p^*\Tang_{\mathcal{X}}$ is equivalent to the $2$-term complex 
    \begin{equation*}
     \hofib\left(T_{X_1|_{X_0}} \overset{t_*-s_*}{\longrightarrow}T_{X_0}\right)  \simeq \big(\ker(s_*)[1] \xrightarrow{t_*} T_{X_0}\big)\,.
    \end{equation*}
    (These are quasi-isomorphic because $s_*$ is assumed to be surjective).    
    Recall that $L=\ker(s_*)$ is the Lie algebroid of the Lie groupoid $X_\bullet$, with anchor map $a=t_{*|\ker(s_*)}$. 
    This $2$-term complex therefore coincides with the underlying complex of the adjoint representation of the Lie algebroid $L$ from \cite{AAC}. 
    Its structure of $L$-representation (up to homotopy) is the shadow of the fact that it is obtained by pulling back along $p:X_0\to\mathcal X$. 
\item Dually, we get that $p^*\Cotang_{\mathcal{X}}$ is equivalent to the $2$-term complex 
\[
\Omega^1_{X_0}\overset{a^*}{\longrightarrow}L^*\,. 
\]
Therefore, 
    \begin{equation*}
        \Sym^2(p^*\Cotang_{\mathcal{X}}[-1]) \simeq 
        \Sym^2\big(\Omega^1_{X_0}[-1] \xrightarrow{a^*}L^*[-2]\big).
    \end{equation*}
\end{enumerate}

Assume we have a closed $2$-form $\omega=\omega_0+\cdots+\omega_n$ of degree $n$ on $\mathcal{X}$. 
We are back to the discussion from \cref{section:linear}, and thus there are only four situations ensuring that the underlying $2$-form is non-degenerate: 
\begin{enumerate}[label=(\alph*)]
    \item If the anchor map is an isomorphism, then $p^*\Cotang_{\mathcal{X}}\simeq 0$. Hence $\Cotang_{\mathcal{X}} \simeq 0$, and there is only the zero form $\omega=0$, which is shifted symplectic for any shift. 
    \item If the anchor $a$ is surjective, then $\g=\ker (a)$ is a bundle of Lie algebras over $X_0$ and $p^*\Cotang_{\mathcal{X}}\simeq \g^*[-1]$ sits in degree $1$ . In this case we have seen that we can only have a degree 2 symplectic form. 
    Consider the subcomplex $(\widehat{DR}^{\geq2}(\mathcal{X}), \widecheck{\partial}+ d)$ of the totalization of the de Rham complex (by definition, cocycles in this complex are closed $2$-forms). Then one can show that the diagram
    \[\begin{tikzcd}
	{\Gamma\left(X_0, \Sym^2(\g^*[-2])\right)^{X_1}} & {\left(\widehat{DR}^{\geq 2}(\mathcal{X}), \widecheck{\partial} + d\right)} \\
	{\Gamma\left(X_0, \Sym^2(\g^*[-2])\right)} & {\Gamma\left(X_0, \Sym^2(p^*\Cotang_{\mathcal{X}}[-2])\right)}
	\arrow[from=1-1, to=1-2]
	\arrow[hook, from=1-1, to=2-1]
	\arrow[from=1-2, to=2-2]
	\arrow["\sim", from=2-1, to=2-2]
    \end{tikzcd}\]
    commutes. Here the top horizontal map is the map 
    \begin{equation*}
        \langle\cdot, \cdot\rangle \longmapsto 
        \frac{1}{2} \langle pr_1^*\theta^L, pr_2^*\theta^R\rangle \in \Omega^2(X_1 \times_{X_0} X_1),
    \end{equation*}
    associating to each $X_1$-invariant pairing on $\g^*$ half of its evaluation on the pullbacks of the Maurer-Cartan forms $\theta^L, \theta^R$. 
    Note that the bottom horizontal map is a quasi-isomorphism. The left vertical map is just the inclusion, while the right vertical map is
    \begin{equation*}
        \omega \longmapsto p^*\omega_0.
    \end{equation*}

    An example of this situation is when $X_0=*$. Then $X_1=G$ is a Lie group and $\mathcal X=[*/G]=BG$. In this case, any invariant metric on $\g$ defines a $2$-shifted symplectic structure on $BG$.
    \item If the anchor is injective with constant rank, then $p^*\Cotang_{\mathcal{X}}\simeq\ker(a^*)$ sits in degree 0 and we can only have a degree $0$ form $\omega_0 \in \Omega^2(X_0)$, satisfying $d\omega =0$ and $\widecheck{\partial} \omega =0$ (in particular, $\omega$ is constant on the leaves of the regular foliation described by the image of $a$). The non-degeneracy condition says that $\omega_0^\flat:\coker(a)\to\ker(a^*)$ is an isomorphism, i.e.~that $\omega_0$ is transversally non-degenerate. 
    \item If none of the above are true, then $\omega$ must have degree 1 and we recover Ping Xu's notion of a quasi-symplectic groupoid \cite{Xu}, where $\omega = \omega_0 + \omega_1$, with $\omega_0\in \Omega^2(X_1)$ and $\omega_1\in \Omega^3(X_0)$. 
\end{enumerate}

The above discussion has the following interesting consequence: in order to have a $0$-shifted symplectic structure on a geometric $1$-stack 
$\mathcal X=[X_0/X_1]$ having non discrete isotropy (for instance in situation (d)), the stack $\mathcal X$ (and in particular $X_0$) must be derived. 
Indeed, the cotangent complex will be non-trivial in degree $-1$, and thus the non-degeneracy 
condition imposes that it will also be non-trivial in degree $1$, which can't occur without a non-trivial derived enhancement of $X_0$. This happens for instance with 
the moduli stack of $G$-local systems on a closed surface (see \cref{example-poincare} and \cref{example-final} below), and is analogous to what is observed in the 
classical BV-BRST formalism where every Chevalley generator (also known as a ghost variable) gets paired with a dual Tate generator (we refer to 
\cite{Anel-Calaque} for more details on this matter). 

\subsection{Lagrangian structures}

The definition of a lagrangian structure for a morphism of stacks is the same as the one for a morphism of derived affine schemes (\cref{definition:lagrangian-affine}). 

\begin{Definition}
Let $f:\mathcal{Y}\to\mathcal{X}$ be a map of (possibly derived) stacks. 
Let $\omega$ be an $n$-shifted symplectic structure on $\mathcal{X}$.     
A \emph{lagrangian structure} on $f$ is a homotopy $\eta$ between $f^*\omega$ and $0$ such that $\eta_0$ is non-degenerate as a \emph{linear} isotropic structure for $f_*: \Tang_{\mathcal{Y}} \to f^*\Tang_{\mathcal{X}}$ with respect to $f^*\omega$.
\end{Definition}

\begin{Example}[Moment maps as lagrangian morphisms]
Let $X_\bullet$ be a Lie groupoid and let $f_\bullet: Y_\bullet \to X_\bullet$ be an action of $X_\bullet$ on $Y_0$, i.e. we assume that $Y_\bullet$ is a Lie groupoid and that 
\begin{equation*}
    Y_n \simeq Y_0 \times_{X_0} X_n,
\end{equation*}
or in other words, $Y_\bullet$ is the action Lie groupoid of the action of $X_\bullet$ on $Y_0$. 
In this case, if $L$ is the Lie algebroid of $X_\bullet$, where we recall that on $X_0$ we have the complex of vector bundles
\begin{equation*}
    p^*\Tang_{|X_\bullet|} \simeq( L[1] \xrightarrow{a} T_{X_0})
\end{equation*}
concentrated in degrees $-1$ and 0, then $f_0^*L$ is the Lie algebroid of $Y_\bullet$. 
This implies that we are in the following situation
\[\begin{tikzcd}
	E & E \\
	B & F
	\arrow[equals, from=1-1, to=1-2]
	\arrow[from=1-1, to=2-1]
	\arrow[from=1-2, to=2-2]
	\arrow[from=2-1, to=2-2]
\end{tikzcd}\]
from \cref{subsec-2.3.4} on the level of tangent complexes.

Assume that we have a quasi-symplectic structure $(\omega_0, \omega_1) \in \Omega^2(X_1) \oplus \Omega^3(X_0)$ on $X_\bullet$. 
Then $\omega_0 + \omega_1$ defines a $1$-shifted symplectic structure on $|X_\bullet|$. We now determine what a lagrangian structure on $|f_\bullet|:|Y_\bullet| \to |X_\bullet|$ is.

First of all, an isotropic structure for $|f_\bullet|$ is a $2$-form $\eta=\eta_0 \in \Omega^2(X_0)$ such that 
\begin{equation*}
    (\widecheck{\partial} + d)(\eta_0) = f_0^*\omega_0 \text{ and } d\eta_0 = f_1^*\omega_1.
\end{equation*}
Which is equivalent to say that $Y_0$ is a pre-hamiltonian space in the sense of \cite[Definition 3.1]{Xu}. 

As we previously mentioned in \cref{subsec-2.3.4}, the non-degeneracy condition coincides with Xu's condition for hamiltonian $X_1$-spaces. \hfill$\triangle$
\end{Example}

The yoga of lagrangian correspondences remains the same. Let us consider some practical examples of this.

\begin{Example}\label{ex:symp-reduction-lagrangian-intersection}
    Let $G$ be a Lie group (or an affine algebraic group). 
    Let $T^*G \cong \g^*\times G \rightrightarrows \g^*$ be the action groupoid of the coadjoint action. This is a symplectic groupoid with respect to the canonical symplectic form on $T^*G$. Therefore the stack $[\g^*/G]$ that it presents is $1$-shifted symplectic. Observe that by the discussion above, hamiltonian $G$-spaces in the classical sense are really hamiltonian $T^*G$-spaces in the sense of \cite{Xu}. Moreover, Marsden--Weinstein reduction can be obtained by a lagrangian intersection similar to the one in \cref{ex:1-shifted-lagrangian-intersection}, as we explain now. 
    Let $\mu: M\to \g^*$ be a momentum map, making $M$ a hamiltonian $G$-space. Then the induced map $[M/G] \to [\g^*/G]$ between quotient stacks carries a lagrangian structure. Because any coadjoint orbit $\mathcal{O}$ is a hamiltonian $G$-space, the map $[\mathcal{O}/G] \to [\g^*/G]$ also carries a lagrangian structure. Thus, by intersecting these two lagrangians we get
    \[\begin{tikzcd}[column sep=scriptsize]
        {M/\!/_{_\mathcal{O}}G} & {[\mathbf{R}\mu^{-1}(\mathcal{O})/G]} & {[\mathcal{O}/G]} \\
        & {[M/G]} & {[\g^*/G],}
        \arrow["{:=}"{description}, draw=none, from=1-1, to=1-2]
        \arrow[from=1-2, to=1-3]
        \arrow[from=1-2, to=2-2]
        \arrow["\lrcorner"{anchor=center, pos=0.125}, draw=none, from=1-2, to=2-3]
        \arrow[from=1-3, to=2-3]
        \arrow[from=2-2, to=2-3]
    \end{tikzcd}\]
    where $\mathbf{R}\mu^{-1}(\mathcal O)$ is the derived fiber of $\mathcal O$ along $\mu$.  Thus $M/\!/_{_\mathcal{O}}G$ carries a $0$-shifted symplectic structure. Whenever $\mathcal O$ is the orbit of a regular value of the moment map, the derived fiber coincides with the underived one $\mu^{-1}(\mathcal{O})$ and we recover the usual Marsden--Weinstein symplectic reduction from \cite{Marsden-Weinstein}. \hfill$\triangle$
\end{Example}

We show how to compute $M/\!/_{_\mathcal{O}}G$ in the above example in the case of $M=\Spec(A)$ a smooth underived affine scheme (this situation is the algebro-geometric analog of $X$ being an honest manifold) and $\mathcal O=\{0\}$. 
In this case we have $\mathbf{R}\mu^{-1}(0) = \Spec(B)$, where $B:= A\overset{\mathbf{L}}{\underset{\Sym(\g)}{\otimes}}k$. Using a quasi-free resolution as in \cref{ex-long}, one gets that $B$ is equivalent to $A \otimes \Sym(\g [1])$ equipped with differential $\delta$ described as follows: 
\begin{equation*}
    \delta(a \otimes 1) = 0, \text{ for } a \in A, \text{ and }
    \delta(1\otimes x) = \mu^*x \otimes 1, \text{ for } x \in \g,
\end{equation*}
where $x$ is seen as a linear function on $\g^*$.
Since all algebras involved are $G$-algebras and all maps are $G$-equivariant, $B$ also carries a $G$-action, given by the action on $A$ tensored with the adjoint action on $\Sym(\g [1])$. 

Then there is a map from the Cartan (graded mixed) complex 
\begin{equation*}
    \left(\Sym_B^\bullet(\Omega^1_B[-1]) \otimes \Sym(\g^*[-2])\right)^G, 
\end{equation*}
with internal differential given by $\delta + d_C$ where $d_C$ is the Cartan differential, and the de Rham differential given by the one on $B$, to $DR\big([\Spec(B)/G]\big)$.

Let $\{x_i\}_i$ be a basis of $\g$. Then the reduced symplectic structure is $\omega_A + \sum_{i} dx_i \cdot x_i^*$. The tangent complex of $M/\!/_{_{\mathcal{O}}} G$ is the $B\rtimes G$-module
\begin{equation*}
    \g \otimes B \rightarrow T_A \otimes_A B \xrightarrow{T_\mu} \g^* \otimes B,
\end{equation*}
in degrees $-1$ to $1$. 

\begin{Example}
    If $X_\bullet$ is a symplectic groupoid, then the trivial isotropic structure on $X_0 \to |X_\bullet|$ is non-degenerate. 
    Thus, in \cref{ex:symp-reduction-lagrangian-intersection}, the quotient map $\g^* \to [\g^*/G]$ is lagrangian. The intersection of this with $[M/G] \to [\g^*/G]$ recovers the $0$-shifted symplectic structure on $M$:
    \[\begin{tikzcd}[column sep=scriptsize]
        M & {\g^*} \\
        {[M/G]} & {[\g^*/G].}
        \arrow[from=1-1, to=1-2]
        \arrow[from=1-1, to=2-1]
        \arrow["\lrcorner"{anchor=center, pos=0.125}, draw=none, from=1-1, to=2-2]
        \arrow[from=1-2, to=2-2]
        \arrow[from=2-1, to=2-2]
    \end{tikzcd}\]
    (see e.g.~\cite{Safronov}). \hfill$\triangle$
\end{Example}

\begin{Example}\label{example:q-hamiltonian-intersection}
    Let $G$ be a Lie group with an invariant metric on $\g$. The action groupoid $G\times G \rightrightarrows G$ with respect to the adjoint action is quasi-symplectic (see \cite[Proposition 2.8]{Xu}). That is, in particular, $[G/G]$ is $1$-shifted symplectic. Observe that, in analogy with \cref{ex:symp-reduction-lagrangian-intersection}, quasi-hamiltonian $G$-spaces in the sense of \cite{AMM} are hamiltonian $(G \times G \rightrightarrows G)$-spaces in the sense of \cite{Xu}. 
    Let $\mu: M\to G$ be a $G$-valued momentum map, giving $M$ the structure of a quasi-hamiltonian $G$-space. Then the induced map $[M/G] \to [G/G]$ on the level of quotient stacks carries a lagrangian structure. Any conjugacy class $\mathcal{C}$ is a quasi-hamiltonian $G$-space as well, which implies $[\mathcal{C}/G] \to [G/G]$ also carries a lagrangian structure. The derived intersection of these two lagrangian is the reduction $M/\!/_{_\mathcal{C}}G$, which is indeed $0$-shifted symplectic:
    \[\begin{tikzcd}[column sep=scriptsize]
        {M/\!/_{_\mathcal{C}}G} & {[\mathbf{R}\mu^{-1}(\mathcal C)/G]} & {[\mathcal{C}/G]} \\
        & {[M/G]} & {[G/G].}
        \arrow["{:=}"{description}, draw=none, from=1-1, to=1-2]
        \arrow[from=1-2, to=1-3]
        \arrow[from=1-2, to=2-2]
        \arrow["\lrcorner"{anchor=center, pos=0.125}, draw=none, from=1-2, to=2-3]
        \arrow[from=1-3, to=2-3]
        \arrow[from=2-2, to=2-3]
    \end{tikzcd}\]
As in \cref{ex:symp-reduction-lagrangian-intersection}, if $\mathcal C$ is the conjugacy class of a regular value of $\mu$, then we recover the quasi-hamiltonian reduction of Alekseev--Malkin--Meinrenken \cite{AMM}. \hfill$\triangle$
\end{Example}
\begin{Remark}
    The quasi-hamiltonian reduction procedure from \cite[Theorem 5.1]{AMM} is actually a bit more general, and can also be recovered using the yoga of compositions of lagrangian correspondences. Let $M$ be a quasi-hamiltonian $G_1\times G_2$-space ($G_1$ and $G_2$ being Lie groups), that induces a lagrangian correspondence 
    \[
    \begin{tikzcd}[sep=scriptsize]
	& {[M/G_1\times G_2]} &\\
	{[G_1/G_1]}& & {[G_2/G_2]}\\
	\arrow[from=1-2, to=2-1]
	\arrow[from=1-2, to=2-3]
\end{tikzcd}
    \]
Given a conjugacy class $\mathcal C\subset G_2$, inducing a lagrangian morphism $[\mathcal{C}/G_2] \to [G_2/G_2]$, we can use the composition of lagrangian correspondences
    
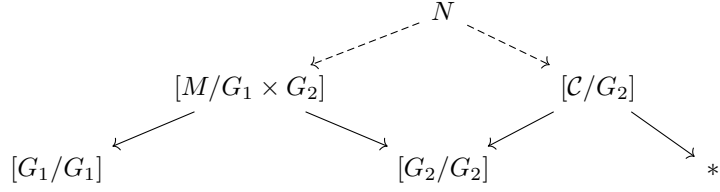
\[\begin{tikzcd}[sep=scriptsize]
	&& {N} \\
	& {[M/G_1\times G_2]} && {[\mathcal{C}/G_2]} \\
	{[G_1/G_1]} && {[G_2/G_2]} && {{*}}
	\arrow[dashed, from=1-3, to=2-2]
	\arrow[dashed, from=1-3, to=2-4]
	\arrow[from=2-2, to=3-1]
	\arrow[from=2-2, to=3-3]
	\arrow[ from=2-4, to=3-3]
	\arrow[from=2-4, to=3-5]
\end{tikzcd}\]
and get a lagrangian morphism $N\to [G_1/G_1]$. One can prove that $N$ is the quotient of $M/\!/_{_\mathcal{C}}G_2=[\mathbf{R}\mu_2^{-1}(\mathcal C)/G_2]$ by $G_1$, where $\mu_2:M\to G_2$ is the $G_2$-valued moment map. In other words, the $G_2$-reduction $M/\!/_{_\mathcal{C}}G_2$ of $M$ is a quasi-hamiltonian $G_1$-space. 
\end{Remark}


\section{AKSZ/PTVV construction}\label{section5}


\subsection{Mapping stacks}

\begin{Definition}
    Let $\mathcal{X}$ and $\mathcal{Y}$ be stacks. The \emph{mapping stack} between $\mathcal{X}$ and $\mathcal{Y}$ is the stack $\Map(\mathcal{X},\mathcal{Y})$, defined on $\Spec A$ by
    \begin{equation*}
        \Map(\mathcal{X},\mathcal{Y})(\Spec A) := \hom_{\mathsf{dSt}}(\mathcal{X} \times \Spec A, \mathcal{Y}),
    \end{equation*}
    where $\hom_{dSt}$ is the hom-space in the ($\infty$-)category of derived stacks. 
\end{Definition}
\begin{Remark}\label{remark:tangent}
 By definition, a map $x: \Spec A \to \Map (\mathcal X,\mathcal Y)$ is the same as a map $f_x : \mathcal X \times \Spec A \to \mathcal Y$. 
 Since derived stacks form an $\infty$-category, there is a space/$\infty$-groupoid of such maps. 
\end{Remark}

\begin{Definition}
    Let $X$ be an $\infty$-groupoid, for instance any CW complex. The \emph{Betti stack} $X_B$ is the sheaf associated to the constant presheaf $X$. 
\end{Definition}

If $X$ is contractible, then $X_B = *$. If $X$ is a CW complex, then $X$ is a gluing (i.e. a colimit) of its cells. All cells are contractible, hence $X$ is a homotopy colimit of a diagram of points. Since the functor $X \mapsto X_B$ preserves homotopy colimits, $X_B$ is also a homotopy colimit (in stacks) of a diagram of points. 

This makes mapping stacks out of Betti stacks easy to calculate, because $\Map(-,\mathcal Y)$ sends homotopy colimits to homotopy limits. Therefore, if $X$ is a CW complex, then $\Map(X_B,\mathcal Y)$ is a homotopy limit (and the diagram is made of products of $\mathcal Y$'s with morphisms being made of diagonals and projections). 

\begin{Example}\label{example:circle}
The circle $S^1$ can be obtained by gluing two $1$-cells along two $0$-cells: 
\begin{equation*}
    S^1 \cong \clim\left(\DecompCircle\right)\simeq\colim\left(\DecompCircle\right) \simeq \colim\left(\BettiCircle\right).
\end{equation*}
Therefore $S^1_B\simeq \colim\left(\BettiCircleb\right)$, and 
\[
\Map(S^1_B,\mathcal Y)\simeq \holim\left(\begin{tikzcd}[sep=tiny]
	{\mathcal Y} \\
	&& \begin{array}{c} \begin{array}{c}{\mathcal Y}\\\times \\{\mathcal Y}\end{array} \end{array} \\
	{\mathcal Y}
	\arrow[from=1-1, to=2-3]
	\arrow[from=3-1, to=2-3]
    \end{tikzcd}\right)\simeq \mathcal Y\underset{\mathcal Y\times\mathcal Y}{\overset{h}{\times}}\mathcal Y.
\]
In the case when $X=BG$, for a Lie (or affine algebraic) group $G$, we recover that $\Map(S^1_B,BG)\simeq[G/G]$. Indeed, the following sequence of homotopy pull-back squares 
\[\begin{tikzcd}
	{G} & {*} \\
	{\Map(S^1_B,BG)} & {[*/G]} \\
	{[*/G]\simeq [G/G\times G]} & {[*/G\times G]}
	\arrow[from=1-1, to=1-2]
	\arrow[from=1-1, to=2-1]
	\arrow["\lrcorner"{anchor=center, pos=0.125}, draw=none, from=1-1, to=2-2]
	\arrow[from=1-2, to=2-2]
	\arrow[from=2-1, to=2-2]
	\arrow[from=2-1, to=3-1]
	\arrow["\lrcorner"{anchor=center, pos=0.125}, draw=none, from=2-1, to=3-2]
	\arrow[from=2-2, to=3-2]
	\arrow[from=3-1, to=3-2]
\end{tikzcd}\]
exhibits $\Map(S^1_B,BG)$ as the quotient of $G$ by an action of $G$ (it is an \textit{exercise} to check that this is the adjoint action). \hfill$\triangle$
\end{Example}

\begin{Example}\label{example:free group}
Let $\Sigma$ be an oriented surface of genus $g$, and let  $\mathring{\Sigma}:=\Sigma-\mathbb{D}^2$ be the surface minus a small disc (in other words, $\Sigma$ is obtained from the $\mathring{\Sigma}$ by ``gluing a $2$-cell'' along the boundary). Here is an example in genus $2$:  
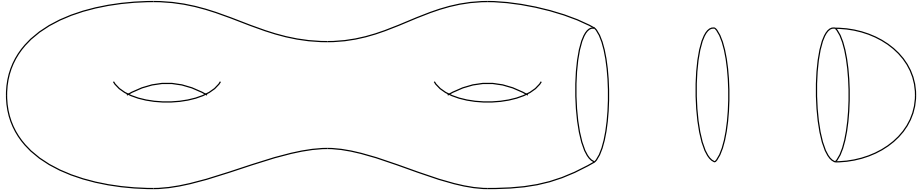
\begin{center}
\begin{adjustbox}{width=\textwidth}
\begin{tikzpicture}[scale=0.6]
	\begin{pgfonlayer}{nodelayer}
		\node [style=none] (1) at (3, 0) {};
		\node [style=none] (2) at (5.75, 1.75) {};
		\node [style=none] (3) at (5.75, -1.75) {};
		\node [style=none] (4) at (9, 1) {};
		\node [style=none] (5) at (12, 1.75) {};
		\node [style=none] (6) at (14, 1.25) {};
		\node [style=none] (7) at (14, -1.25) {};
		\node [style=none] (8) at (9, -1) {};
		\node [style=none] (9) at (12, -1.75) {};
		\node [style=none] (10) at (5, 0.25) {};
		\node [style=none] (11) at (7, 0.25) {};
		\node [style=none] (12) at (5.25, 0) {};
		\node [style=none] (13) at (6.75, 0) {};
		\node [style=none] (14) at (11, 0.25) {};
		\node [style=none] (15) at (13, 0.25) {};
		\node [style=none] (16) at (11.25, 0) {};
		\node [style=none] (17) at (12.75, 0) {};
		\node [style=none] (18) at (16.25, 1.25) {};
		\node [style=none] (19) at (16.25, -1.25) {};
		\node [style=none] (20) at (18.5, 1.25) {};
		\node [style=none] (21) at (18.5, -1.25) {};
		\node [style=none] (22) at (20, 0) {};
	\end{pgfonlayer}
	\begin{pgfonlayer}{edgelayer}
		\draw [in=180, out=-90] (1.center) to (3.center);
		\draw [in=-180, out=90] (1.center) to (2.center);
		\draw [in=-180, out=0, looseness=0.75] (3.center) to (8.center);
		\draw [in=180, out=0] (2.center) to (4.center);
		\draw [in=180, out=0] (4.center) to (5.center);
		\draw [in=150, out=0, looseness=0.75] (5.center) to (6.center);
		\draw [in=180, out=0, looseness=0.75] (8.center) to (9.center);
		\draw [in=-150, out=0] (9.center) to (7.center);
		\draw [in=165, out=165, looseness=0.50] (6.center) to (7.center);
		\draw [in=45, out=-45, looseness=0.50] (6.center) to (7.center);
		\draw [bend right=60, looseness=0.75] (10.center) to (11.center);
		\draw [bend left] (12.center) to (13.center);
		\draw [bend right=60, looseness=0.75] (14.center) to (15.center);
		\draw [bend left] (16.center) to (17.center);
		\draw [in=165, out=165, looseness=0.50] (18.center) to (19.center);
		\draw [in=45, out=-45, looseness=0.50] (18.center) to (19.center);
		\draw [in=165, out=165, looseness=0.50] (20.center) to (21.center);
		\draw [in=45, out=-45, looseness=0.50] (20.center) to (21.center);
		\draw [in=90, out=0] (20.center) to (22.center);
		\draw [in=-90, out=0] (21.center) to (22.center);
	\end{pgfonlayer}
\end{tikzpicture}
\end{adjustbox}
\end{center}
We get that 
$\Sigma \simeq \mathring{\Sigma}\underset{S^1}{\sqcup} \mathbb{D}^2\simeq \mathring{\Sigma}\overset{h}{\underset{S^1}{\sqcup}} *$. 
Now recall that $\mathring{\Sigma}$ is weakly equivalent to a wedge of $2g$ circles, and thus $\mathring{\Sigma}_B$ can be identified with $BF_{2g}$, where $F_{2g}$ is the free group on $2g$ generators. Moreover, for an appropriate choice of generators $a_1,\dots,a_g,b_1,\dots,b_g$ giving the identification $F_{2g}\cong\pi_1(\mathring{\Sigma})$, the map $B\mathbb{Z}\simeq S^1_B\to \mathring{\Sigma}_B\simeq BF_{2g}$ is induced by the group morphism 
\begin{eqnarray*}
    \mathbb{Z} & \longrightarrow & F_{2g} \\
    1 & \longmapsto & \vec{\prod_{i=1,\dots g}}[a_i,b_i].
\end{eqnarray*} 
Hence we have 
\begin{equation*}
\begin{split}
    \Map(\Sigma_B, BG) 
    &\simeq \Map\left(\mathring{\Sigma}_B \overset{h}{\underset{S^1_B}{\sqcup}} *\,, BG\right) 
    \simeq \Map(\mathring{\Sigma}_B) \overset{h}{\underset{[G/G]}{\times}} [*/G] \\
    & \simeq [G^{2g}/G] \overset{h}{\underset{[G/G]}{\times}} [*/G] \simeq [\mathbf{R}\mu^{-1}(1)/G],
\end{split}
\end{equation*}
where $\mu:G^{2g}\to G$ sends $(A_1,\dots,A_g,B_1,\dots,B_g)$ to the product of commutators $\prod_i(A_i,B_i)$. This recovers (a derived enhancement of) the description of the character variety/stack of a closed oriented surface as a quasi-hamiltonian reduction from \cite{AMM} (see also \cite{Meinrenken}). 
Note that in the case of the $2$-sphere, $\mathbf{R}\mu^{-1}(1)\simeq\mathfrak{g}[-1]$, so that $\Map(S^2_B, BG)\simeq\big[\mathfrak{g}[-1]/G\big]$ while the 
\textit{underived} moduli stack of $G$-local systems on the $2$-sphere is $[*/G]$. 
\hfill$\triangle$
\end{Example}

In the remaining sections, we explain how the above example still fits well within the shifted symplectic geometry approach. 


\subsection{Shifted symplectic structures on mapping stacks}

\begin{Definition}\label{definition:orientation}
    If $\mathcal X$ is a stack, then $[\mathcal X]: \mathbf{R}\Gamma(\mathcal X, \mathcal{O}_{\mathcal X}) \to k[-d]$ is a $d$-orientation if:
    \begin{enumerate}[label=(\alph*)]
        \item For any perfect complex $E$ on $\mathcal X \times \Spec A$, the (derived) sections space $\mathbf{R}\Gamma(\mathcal X \times \Spec A, E)$ is a perfect $A$-module.
        \item For any perfect complex $E$ on $\mathcal X$, the map
        \begin{equation*}
        \begin{split}
            \mathbf{R}\Gamma(\mathcal X, E^*) &\longrightarrow \mathbf{R}\Gamma(\mathcal X, E)^*[-d]\\
            \xi &\longmapsto (s \mapsto [X](s(\xi)))
        \end{split}
        \end{equation*}
        is a quasi-isomorphism. 
    \end{enumerate}
\end{Definition}

If $X$ is a space, then sheaves on $X_B$ are local systems of complexes on $X$, and for such a local system $E$, there is an equivalence $\mathbf{R}\Gamma(X_B,E)\simeq C_{sing}^\bullet(X,E)$. Hence if $X$ has finite homotopy type, then $X_B$ satisfies condition (a) of the definition above. Condition (b) for $X_B$ is then equivalent to requiring that $X$ is a Poincaré duality space of dimension $d$ (the $d$-orientation is given by using the cap-product with the fundamental class of $X$). 

\begin{Example}\label{example:circle2}
If $X$ is a CW complex, then $\mathbf{R}\Gamma(X_B, \mathcal{O}_{X_B})\simeq C_{cell}^\bullet(X, k)$ is the complex of cellular cochains of $X$. For example, using the CW decomposition of $S^1_B$ from \cref{example:circle} we get: 
\begin{equation*}
    \mathbf{R}\Gamma(S^1_B,\mathcal{O}_{S^1_B}) \simeq  \hofib\left(
    \begin{array}{lll}
        k\oplus k &\to &k\oplus k \\
        (x,y) &\mapsto &(x-y, y-x)
    \end{array}
    \right).
\end{equation*}
Since the fundamental class of $S^1$ (recall that $S^1$ is a Poincaré duality space of dimension $1$) is given by the sum of the two $1$-cells in the cellular decomposition from \cref{example:circle}, the morphism of $2$-term complexes 
\begin{equation*}
    \left(
    \begin{tikzcd}
    k\oplus k  \\
    k\oplus k
    \arrow[from=1-1, to=2-1]
    \end{tikzcd}
    \right)
    \xrightarrow{\begin{pmatrix}0\\+\end{pmatrix}}    
    \left(\begin{tikzcd}
	0 \\
	k
	\arrow[from=1-1, to=2-1]
    \end{tikzcd}\right)=k[-1] 
    \end{equation*}
provides a description of the induced $1$-orientation on $S^1_B$. \hfill$\triangle$\end{Example}

\begin{Theorem}[{Pantev--Toen--Vaquié--Vezzosi \cite[Theorem 2.5]{PTVV}}]\label{theorem:ptvv}
    If $\mathcal X$ has a $d$-orientation and $\mathcal Y$ has a non-degenerate closed $n$-shifted $2$-form $\omega: \wedge^2\Tang_{\mathcal Y} \to \mathcal{O}_{\mathcal Y}[n]$, then there is a non-degenerate closed $(n-d)$-shifted $2$-form 
    \begin{equation*}
        \int_{[\mathcal X]}\omega_0
    \end{equation*}
    on $\Map(\mathcal X, \mathcal Y)$. 
\end{Theorem}

We refer to \cite{PTVV} for the proof, but we describe here what the underlying $2$-form $(\int_{[\mathcal X]}\omega)_0$ looks like, which actually only depends on the $d$-orientation $[\mathcal X]$ and the underlying $2$-form $\omega_0$. 
First, the following expected result computes the tangent complex of $\Map(\mathcal X,\mathcal Y)$: 
\begin{Lemma}[{see e.g.~\cite[Proposition B.10.21]{CHS}}]
    If $\mathcal X$ satisfies property (a) of \cref{definition:orientation} and $\mathcal Y$ has a perfect (co)tangent complex, then for every $A$-point $x: \Spec A \to \Map(\mathcal X,\mathcal Y)$  given by $f_x : \mathcal X \times \Spec A \to \mathcal Y$ (see \cref{remark:tangent}), 
            $x^* \Tang _{\Map(\mathcal X,\mathcal Y)} \simeq \mathbf{R}\Gamma(\mathcal X \times \Spec A, f_x^*\Tang_{\mathcal Y})$. 
\end{Lemma}
Then one defines the underlying $2$-form of $x^*\int_{[\mathcal X]}\omega$ as the following pairing: 
\begin{equation*}
\begin{split}
    \wedge^2_A\big(\mathbf{R}\Gamma(\mathcal X \times \Spec A, f^*_x\Tang_Y)\big)
    \xrightarrow{f_x^*\omega_0}
    &\mathbf{R}\Gamma(\mathcal X \times \Spec A, \mathcal{O}_{\mathcal X \times \Spec A})[-n]
    \\
    &\quad\simeq \mathbf{R}\Gamma(\mathcal X, \mathcal{O}_{\mathcal X}) \otimes A[n] 
    \xrightarrow{[\mathcal X] \otimes id_A}
    A[n-d].
\end{split}
\end{equation*}
This describes $(\int_{[\mathcal X]}\omega)_0$ on points. 
Its non-degeneracy can be checked on points, and follows from the non-degeneracy of $\omega_0$ together with 
condition (b) from \cref{definition:orientation}. 

\begin{Example}\label{example-poincare}
For any Poincaré duality space $X$ of dimension $d$, we have seen that $X_B$ carries a $d$-orientation. 

Therefore, if $G$ is an affine algebraic group together 
with an invariant metric on its Lie algebra $\mathfrak{g}$ (recall that this defines a $2$-shifted symplectic structure on $BG$), then the derived stack $\Map(X_B,BG)$ 
of $G$-local systems on $X$ is $(2-d)$-shifted symplectic. At a $k$-point $x$ of this stack (that is given by a $G$-local system $L_x$), the induced linear $(2-d)$-shifted 
symplectic structure coincides with the one from \cref{ex2} on $x^*\mathbb{T}_{\Map(X_B, BG)}\simeq H^\bullet(X,\ad L_x)[1]$. 

In contrast with the 
\textit{underived} moduli stack of $G$-local systems on $X$, that only depends on the fundamental groupoid of $X$, the \textit{derived} stack 
$\Map(X_B,BG)$ depends on the homotopy type of $X$. For instance, the positive part of its tangent complex (that is induced by the derived structure) at a $k$-point 
$x$ as above is given by the cohomology groups $H^i(X,\ad L_x)$ for $i\geq2$. 

In dimension $1$, we recover that $[G/G]\simeq\Map(S^1_B,BG)$ is $1$-shifted symplectic (as $1=2-1$). We refer to \cite{Safronov} for the proof that this $1$-shifted symplectic structure coincides with the one coming from the quasi-symplectic groupoid $G\times G \rightrightarrows G$ from \cite[Proposition 2.8]{Xu}. 

In dimension $2$, if $\Sigma$ is a closed oriented surface then $\Sigma_B$ carries a $2$-orientation. Thus, if $G$ is as above, then 
$\Map(\Sigma_B, BG)$ is $0$-shifted symplectic. At a smooth $k$-point $x$, $H^\bullet(\Sigma,\ad L_x)=H^1(\Sigma,\ad L_x)$ and the symplectic form thus 
coincides with the one introduced by Atiyah--Bott \cite{AB} (or rather an algebro-geometric analog of it). 
Outside of the smooth locus, the degree $-1$ part $H^0(\Sigma,\ad L_x)$ of the tangent complex 
(that is tangent to the isotropy) is in duality with its degree $1$ part $H^2(\Sigma,\ad L_x)$ (induced by the derived structure). 
The example of the $2$-sphere is very telling: $\Map(S^2_B, BG)\simeq\big[\mathfrak{g}[-1]/G\big]$ (see \cref{example:free group}) only has one $k$-point, 
that is not smooth, and the tangent complex at this point is $H^\bullet(S^2,\mathfrak{g})[1]\simeq \mathfrak{g}[1]\oplus \mathfrak{g}[-1]$ (the degree $0$ pairing is the metric on $\mathfrak{g}$). 
\hfill$\triangle$\end{Example}

We have seen that $\Map(\Sigma_B, BG)$ can be obtained as the quotient by $G$ of the derived fiber at $1$ of the $G$-valued moment map $G^{2g}\to G$ that appears in the quasi-hamiltonian approach of Alekseev--Malkin--Meinrenken (see \cref{example:free group} above). We will see in the remaining section below how the two pictures fit together. 

\begin{Example}
What is interesting with the global approach/picture using derived geometry and mapping stacks is that it is totally intrinsic. Therefore, different ways of presenting the same stack lead to different computations that automatically give equivalent results. 
For instance, instead of viewing the circle as a gluing of two segments along two points as in \cref{example:circle}, we can also view it as a triangle, that is a gluing of three segments along three points: we glue $x=[a,b]$, $y=[b,c]$ and $z=[c,a]$ along $a$, $b$, and $c$. This amounts to see $S^1_B$ as the classifying stack $BF$ of the groupoid $F$ having objects $\{a,b,c\}$ and freely generating arrows $\{x,y,z\}$. Then, one can prove that for any affine algebraic group $G$ we get  
\begin{equation*}
\Map(S^1_B,BG)\simeq [G^{\{x,y,z\}}/G^{\{a,b,c\}}],
\end{equation*}
were $G=G^{\{i\}}$ (resp.~$G=G^{\{j\}}$) acts by left (resp.~right) multiplication on $G=G^{\{[i,j]\}}$. 
This allows to understand the $1$-shifted symplectic structure on $\Map(S^1_B,BG)$ as a quasi-symplectic groupoid structure on the action groupoid (for the $G^{\{a,b,c\}}$-action on $G^{\{x,y,z\}}$). Since it induces \emph{by construction} the same $1$-shifted symplectic structure on $[G/G]$, it automatically is Morita equivalent to the Alekseev--Malkin--Meinrenken--Xu quasi-symplectic groupoid (for the adjoint action). It is expected that one can get formulas similar to the ones appearing in \cite{Meinrenken} for this quasi-symplectic structure. 

Note that the cochain complex associated with the above CW decomposition is $\hofib(k^{\{a,b,c\}}\to k^{\{x,y,z\}})$ where the map sends $\delta_a$ to $\delta_z-\delta_x$, $\delta_b$ to $\delta_x-\delta_y$, and $\delta_c$ to $\delta_y-\delta_z$. 
The $1$-orientation is given by the map to $k[-1]$ sending the three degree $1$ basis elements $\delta_x$, $\delta_y$ and $\delta_z$ to $1$. 

Finally observe that the same reasoning carries over for the CW decomposition of $S^1$ given as the boundary of an $n$-gon. \hfill$\triangle$
\end{Example}

\subsection{Relative version and compatibility with gluings}

There is a ``relative'' version of the above theorem of Pantev--Toën--Vaquié--Vezzosi, for which we refer to \cite[Theorem 2.11]{Cal}. In the case of Betti stacks, it specializes to the following: 
\begin{Theorem}
If $X$ is an oriented manifold of dimension $d+1$, with boundary $\partial X$, and if $\mathcal Y$ is an $n$-shifted symplectic stack, then the restriction morphism 
\[
\Map(X_B,\mathcal Y)\longrightarrow \Map\big((\partial X)_B,\mathcal Y\big)
\]
carries a lagrangian structure with corresponding 
$(n-d)$-shifted symplectic structure on $\Map\big((\partial X)_B,\mathcal Y\big)$ being given by \cref{theorem:ptvv}. 
\end{Theorem}
The more general version \cite[Theorem 2.11]{Cal} makes use of a relative version of the notion of $d$-orientations, that we will not explain in these notes. 
\begin{Example}
    Let $\mathring{\Sigma}$ and $G$ be as in \cref{example:free group}, and assume that there is a $G$-invariant metric on the Lie algebra $\mathfrak{g}$ (inducing a $2$-shifted symplectic structure on $BG$) Then the morphism 
    \[
    [G^{2g}/G]\simeq\Map(\mathring{\Sigma}_B,BG)\to \Map(S^1_B,BG)\simeq [G/G]
    \]
    carries a lagrangian structure. 
    This lagrangian structure happens to coincide with the one induced by the quasi-hamiltonian $G$-space structure on $G^{2g}$ from \cref{example:q-hamiltonian-intersection}. \hfill$\triangle$
\end{Example}
One way of proving the final claim in the above example, and that \[\Map(\Sigma_B,BG)\simeq G^{2g}/\!/_{_{\{1\}}} G\] \textit{as $0$-shifted symplectic stacks} (recall that for the underlying stacks, i.e.~without the symplectic structure, this follows from \cref{example:free group}), is to prove that both the mapping stack construction and the quasi-hamiltonian one obey the same ``cut-and-paste properties'' and that they agree on building blocks. 

For the quasi-hamiltonian picture, this follows from \cite{AMM} (see also \cite{Meinrenken}). The fact that they agree on building blocks has been proven in \cite{Safronov}. Finally, the ``cut-and-paste'' properties of the mapping stack construction is summarized in the following result, that we state after defining the necessary categories:
\begin{itemize}
\item Let $\mathsf{Cob}^{or}_d$ be the category with objects closed oriented $(d-1)$-dimensional manifolds, and whose hom spaces are classifying oriented cobordisms between them. This is a monoidal category with respect to the disjoint union $\sqcup$. 
\item Let $\mathsf{Lag}_{n-d+1}$ be the category with objects $(n-d+1)$-shifted symplectic stacks, and hom spaces classifying lagrangian correspondences. This is a monoidal category with respect to the cartesian product $\times$. 
\end{itemize}

\begin{Theorem}
    Every $n$-shifted symplectic stack $\mathcal Y$ defines a symmetric monoidal functor 
    \begin{equation*}
        \Map\big((-)_B, \mathcal Y\big): \mathsf{Cob}^{or}_d \longrightarrow \mathsf{Lag}_{n-d+1}. 
    \end{equation*}
\end{Theorem}

At the level of homotopy categories, this is \cite[Theorem 4.8]{Cal}. At the level of $\infty$-categories this follows from the main results of \cite{CHS} (more precisely Theorem C and Theorem E). 
\begin{Remark}
    The advantage of the $\infty$-categorical statement is that it gives for free an action of diffeomorphism groups on mapping spaces. 
\end{Remark}
\begin{Remark}
    Note that \cite{CHS} proves an even more general statement involving symmetric monoidal $(\infty,n)$-categories, that allows to consider gluing patterns of codimension higher than $1$ (or, iterated gluing patterns). 
\end{Remark}

\begin{Example}\label{example-final}
For $\mathcal Y=BG$ (with, as usual, $G$ being an affine algebraic group together with an invariant metric on its Lie algebra $\mathfrak g$), $n=2$, $d=2$, we can compute what this functor associates to a genus $g$ closed oriented surface $\Sigma$ by decomposing it as in \cref{example:free group}, interpreted as the composition of $\mathring{\Sigma}$ (viewed as a cobordism from $\emptyset$ to $S^1$) with $\mathbb{D}^2$ (viewed as a cobordism from $S^1$ to $\emptyset$). For instance, with $g=2$ this gives
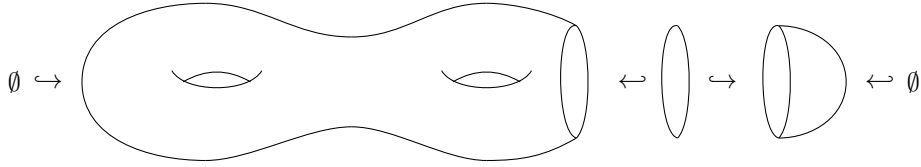
\begin{center}
\begin{adjustbox}{width=\textwidth}
\begin{tikzpicture}[scale=0.6]
	\begin{pgfonlayer}{nodelayer}
		\node [style=none] (0) at (1.5, 0) {$\emptyset$};
		\node [style=none] (1) at (3, 0) {};
		\node [style=none] (2) at (5.75, 1.75) {};
		\node [style=none] (3) at (5.75, -1.75) {};
		\node [style=none] (4) at (9, 1) {};
		\node [style=none] (5) at (12, 1.75) {};
		\node [style=none] (6) at (14, 1.25) {};
		\node [style=none] (7) at (14, -1.25) {};
		\node [style=none] (8) at (9, -1) {};
		\node [style=none] (9) at (12, -1.75) {};
		\node [style=none] (10) at (5, 0.25) {};
		\node [style=none] (11) at (7, 0.25) {};
		\node [style=none] (12) at (5.25, 0) {};
		\node [style=none] (13) at (6.75, 0) {};
		\node [style=none] (14) at (11, 0.25) {};
		\node [style=none] (15) at (13, 0.25) {};
		\node [style=none] (16) at (11.25, 0) {};
		\node [style=none] (17) at (12.75, 0) {};
		\node [style=none] (18) at (16.25, 1.25) {};
		\node [style=none] (19) at (16.25, -1.25) {};
		\node [style=none] (20) at (18.5, 1.25) {};
		\node [style=none] (21) at (18.5, -1.25) {};
		\node [style=none] (22) at (20, 0) {};
		\node [style=none] (23) at (21.5, 0) {$\emptyset$};
		\node [style=none] (24) at (2.25, 0) {$\hookrightarrow$};
        \node [style=none] (25) at (17.25, 0) {$\hookrightarrow$};
        \node [style=none] (26) at (15.25, 0) {$\hookleftarrow$};
        \node [style=none] (27) at (20.75, 0) {$\hookleftarrow$};
	\end{pgfonlayer}
	\begin{pgfonlayer}{edgelayer}
		\draw [in=180, out=-90] (1.center) to (3.center);
		\draw [in=-180, out=90] (1.center) to (2.center);
		\draw [in=-180, out=0, looseness=0.75] (3.center) to (8.center);
		\draw [in=180, out=0] (2.center) to (4.center);
		\draw [in=180, out=0] (4.center) to (5.center);
		\draw [in=150, out=0, looseness=0.75] (5.center) to (6.center);
		\draw [in=180, out=0, looseness=0.75] (8.center) to (9.center);
		\draw [in=-150, out=0] (9.center) to (7.center);
		\draw [in=165, out=165, looseness=0.50] (6.center) to (7.center);
		\draw [in=45, out=-45, looseness=0.50] (6.center) to (7.center);
		\draw [bend right=60, looseness=0.75] (10.center) to (11.center);
		\draw [bend left] (12.center) to (13.center);
		\draw [bend right=60, looseness=0.75] (14.center) to (15.center);
		\draw [bend left] (16.center) to (17.center);
		\draw [in=165, out=165, looseness=0.50] (18.center) to (19.center);
		\draw [in=45, out=-45, looseness=0.50] (18.center) to (19.center);
		\draw [in=165, out=165, looseness=0.50] (20.center) to (21.center);
		\draw [in=45, out=-45, looseness=0.50] (20.center) to (21.center);
		\draw [in=90, out=0] (20.center) to (22.center);
		\draw [in=-90, out=0] (21.center) to (22.center);
	\end{pgfonlayer}
\end{tikzpicture}
\end{adjustbox}
\end{center}
Applying the functor $\Map\big((-)_B, BG\big)$ we get that the result is the lagrangian intersection
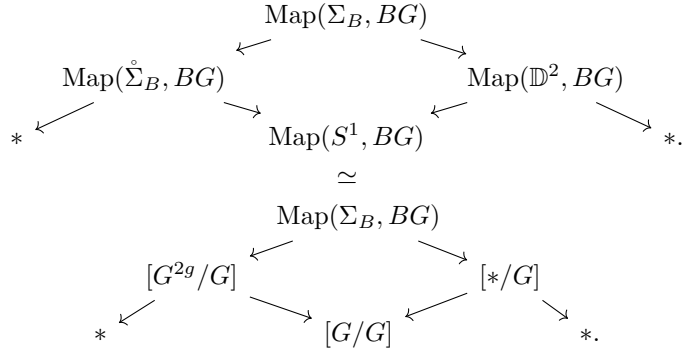
\begin{eqnarray*}
&\begin{tikzcd}[sep=tiny]
	&& {\Map(\Sigma_B, BG)} \\
	& {\Map(\mathring{\Sigma}_B, BG)} && {\Map(\mathbb{D}^2,BG)} \\
	{*} && {\Map(S^1,BG)} && {*.}
	\arrow[from=1-3, to=2-2]
	\arrow[from=1-3, to=2-4]
	\arrow[from=2-2, to=3-1]
	\arrow[from=2-2, to=3-3]
	\arrow[from=2-4, to=3-3]
	\arrow[from=2-4, to=3-5]
\end{tikzcd} &\\
& \simeq & \\
& \begin{tikzcd}[sep=tiny]
	&& {\Map(\Sigma_B, BG)} \\
	& {[G^{2g}/G]} && {[*/G]} \\
	{*} && {[G/G]} && {*.}
	\arrow[from=1-3, to=2-2]
	\arrow[from=1-3, to=2-4]
	\arrow[from=2-2, to=3-1]
	\arrow[from=2-2, to=3-3]
	\arrow[from=2-4, to=3-3]
	\arrow[from=2-4, to=3-5]
\end{tikzcd}&
\end{eqnarray*}
This indeed shows that $\Map(\Sigma_B, BG)\simeq G^{2g}/\!/_{_{\{1\}}}G$ \emph{as $0$-shifted symplectic stacks}. \hfill$\triangle$
\end{Example}

\end{document}